\documentclass[12pt]{report}
\usepackage{amsfonts,amssymb,amsmath,amsbsy}
\usepackage{latexsym}
\usepackage{graphicx}
\usepackage{tikz}

\newtheorem{theo}{Theorem}
\newtheorem{lemma}[theo]{Lemma}
\newtheorem{prop}[theo]{Proposition}

\newtheorem{coro}[theo]{Corollary}
\newtheorem{defi}[theo]{Definition}

\newtheorem{note}[theo]{Note}

\newcommand{\proof}{\noindent{\it Proof: }}
\newcommand{\proofbox}{\hfill \mbox{ $\Box$}\\}

\newcommand{\R}{{\mathbb R}}

\newcommand{\N}{{\mathbb N}}

\newcommand{\M}{{\mathcal{M}}}
\newcommand{\sN}{{\mathcal{N}}}
\newcommand{\ev}{{\rm{Ev}}}
\newcommand{\beg}{{\rm{Beg}}}
\newcommand{\nd}{{\rm{End}}}
\newcommand{\sm}{{\rm{Sim}}}
\newcommand{\be}{{\rm{Bw}}}
\newcommand{\ed}{{\rm{Eq}}}

\newcommand{\cop}{{\rm{Cop}}}

\newcommand{\ftl}{{\rm{FTL}}}
\newcommand{\stl}{{\rm{STL}}}
\newcommand{\lightspeed}{{\rm{Lightspeed}}}
\newcommand{\pow}{\mathcal{P}}
\newcommand{\sftl}{$^{FTL}$}

\title{A logical treatment of special relativity, with and without faster-than-light observers}
\author{Benjamin Hoffman \\ \\ \\ An honors thesis written in partial fulfillment \\ of the requirements for the degree of Bachelors of Arts \\ \\ Department of Mathematical Sciences \\ \\ Lewis \& Clark College \\ \\ Portland, Oregon}

\begin{document}
\maketitle

\renewcommand{\abstractname}{Acknowledgements}
\begin{abstract}
I would like to thank Istv\'an N\'emeti and Hajnal Andr\'eka for suggesting the subject of this thesis, Gergely Sz\'ekely for his extensive comments on drafts of this document, and the entire Algebraic Logic group at the R\'enyi Institute in Budapest for a number of long conversations about logic and relativity. This document would not have been possible were it not for the awesome support of Paul Allen, my ``official advisor''. Biggest thanks go to my ``unofficial advisor'' Iva Stavrov, who pushed me to write something long and mathy, who learned logic with me, and who has spent far longer than is reasonable trying to decipher what I write.
\end{abstract}

\tableofcontents
\chapter{Introduction}

Since Einstein proposed the theory of special relativity in \cite{Einstein}, there have been a number of attempts to formally state its assumptions in a logical axiom system.  Perhaps the most famous of these attempts is the one by Reichenbach (see \cite{Reichenbach}); recently a system has been developed extensively by Andr\'eka, Madar\'asz, Sz\'ekely, and N\'emeti (see \cite{Renyi}, \cite{Renyi2}, and \cite{Renyi3}). Well-developed formal treatments of different aspects of relativity have also been put forth by Goldblatt in \cite{Goldblatt} and Benda in \cite{Benda}.

There are several reasons that this project is interesting. First, by translating special relativity into a formal axiomatization, one removes all ambiguity from its statements, and makes explicit all its assumptions. One may be able to show the relative simplicity of these assumptions (see \cite{Renyi3}), or show that one only needs a handful of undefined predicates to describe many of the features of special relativity (see \cite{Ax}). Second, consequences of special relativity could be proved as theorems that follow from this axiomatization. One could prove novel results in special relativity, potentially using very different methods than those that are usually employed when studying relativity. Alternatively, one could develop the metatheory of relativity. That is, one could prove theorems \emph{about} the theory of special relativity. For instance, one might prove the logical independence of a statement from the theory of special relativity (see \cite{Renyi4} and \cite{Szekely}).

Each of the attempts to axiomatize relativity has different goals, and each comes with its own strengths and weaknesses. For instance, in this thesis we develop an axiom system based on one by Ax in \cite{Ax}. This system is written in a language with only two undefined symbols, and its axioms were designed to be empirically testable and confirmed by a strong body of evidence. Also, unlike the axiom system in the work of Andr\'eka et. al., it does not include a field of quantities to describe the physical objects in its universe. It ontologically commits to the existence of only physical objects: observers and signals. 


It is generally held that the following is an important consequence of special relativity: given any two observers $a$ and $b$, $a$ and $b$'s relative velocities are slower than the speed of light. Recently, in \cite{HillCox} Hill and Cox presented a set of formulas meant to serve as reference frame transformations between two observers who have relative velocities exceeding the speed of light. Hill and Cox cite several motivations for this project, including the possibility of eliminating the use of imaginary masses and complex speeds in some recent physical theories. In \cite{HillCoxNote}, the authors present a simple description of the Hill-Cox transformations, in terms of Poincar\'e transformations. They also prove that the Hill-Cox transformations are consistent with Einstein's principal of relativity only if space is one-dimensional. If we are to admit faster than light observers into our ontology, then, we should examine exactly which of the assumptions we make in special relativity preclude the existence of faster than light observers. These assumptions need to be revised or abandoned.

There are three goals of this thesis. First: to present a concise yet accessible description of basic mathematical logic and model theory. I attempt this in Chapter \ref{LogicChapter}. Second: to develop an axiomatization of special relativity using only two undefined predicates. Ideally, these axioms should be empirically verifiable and supported by a large body of evidence. I give this axiomatization, as well as a complete characterization of its models, in Chapter \ref{STLChapter}. Finally: to weaken this axiomatization so as to allow for the existence of faster-than-light observers, which are normally excluded from the theory of special relativity. The modifications to the axiom system should be as slight as possible, and should be motivated by the question: what features of the original axiom system preclude the existence of faster-than-light observers? I do this in Chapter \ref{FTLchapter}. Throughout this thesis, figures resemble \emph{Minkowski diagrams}, which are the standard way of representing worldlines through spacetime. I give a short introduction to special relativity and Minkowski diagrams in the next chapter.

\chapter{A brief introduction to special relativity and Minkowski diagrams}

In this chapter I present some of the basic features of special relativity, the understanding of which is sufficient to understand the physical basis of this thesis. Any reader interested in a more comprehensive introduction should consult a textbook such as \cite{Schutz}.

In special relativity, spacetime is represented as $\R^4$. Points $\langle x_0,x_1,x_2,x_3\rangle$ in $\R^4$ represent points in spacetime. The first coordinate of such a point gives its position in time, and the remaining three coordinates represent its position in space. A path in $\R^4$ represents the \emph{world line} of some object traveling through spacetime; worldlines run through all the spacetime points that an object moves through in its journey through the universe. Our units are fixed so that any straight world line whose slope has magnitude $1$ (i.e., given points $\langle x_0,x_1,x_2,x_3\rangle,\langle y_0,y_1,y_2,y_3\rangle$ on the line, $\big|\frac{y_0-x_0}{\sqrt{(y_1-x_1)^2+(y_2-x_2)^2+(y_3-x_3)^2}}\big|=1$), represents an object that is moving at the speed of light.

\begin{figure}
\begin{center}
\begin{tikzpicture}[scale=1]
\draw [dashed] (-2,0) -- (2,0) node[right]{space};
\draw [dashed] (0,-2.5) -- (0,2.5) node[right]{time};
\draw (0,-2) node[left]{$a$} -- (0,0) node[below right]{$o$} -- (0,2);
\draw (-.5,-2) node[left]{$b$} -- (-1.75,2);
\draw (.25,-2) node[right]{$c$} -- (1.25,2);
\draw (.75,-2) node[right]{$d$} -- (1.75,2);
\end{tikzpicture}
\end{center}
\caption{A two-dimensional Minkowski diagram of four inertial observers. The observers $c$ and $d$ are moving with the same relative velocity. The diagram represents $a$'s reference frame, since $a$ passes vertically through the origin $o$.} \label{Min1}
\end{figure}
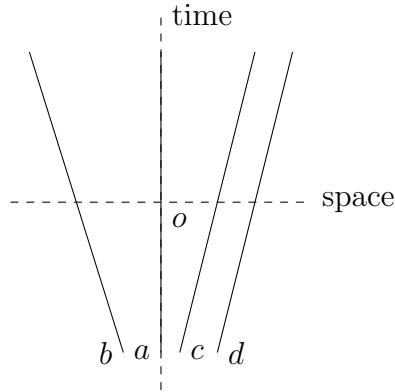

The standard way to represent the geometry of spacetime is using Minkowski diagrams (Figure \ref{Min1}). In a Minkowski diagram, the four coordinate axes (corresponding to the three spatial dimensions, plus time) are drawn such that the time coordinate axis runs vertically. Often, all but one or two of the space axes are not drawn. Objects moving through spacetime are represented by their world lines . An inertial (non-accelerating) object will have a straight line for its world line; two objects moving at the same relative velocity will have parallel world lines. In the remainder of this paper, we restrict discussion to only inertial objects.

We define the \emph{Minkowskian inner product} on $\R^4$. Given vectors $x=\langle x_0,x_1,x_2,x_3\rangle$ and $y=\langle y_0,y_1,y_2,y_3\rangle$, where $x_0$ and $y_0$ represent the time coordinate, the inner product is given:
\begin{equation}
\langle x,y \rangle = -x_0y_0+x_1y_1+x_2y_2+x_3y_3.
\end{equation}
The Minkowskian inner product gives us a way to measure the ``distance'' between two points in spacetime (though this distance can be negative, and two non-identical points may be of distance $0$ from each other). It also gives us a definition of orthogonality. 

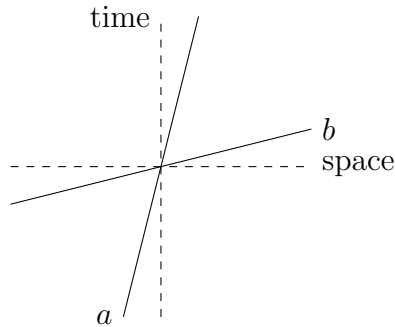
\begin{figure}
\begin{center}
\begin{tikzpicture}[scale=1]
\draw [dashed] (-2,0) -- (2,0) node[right]{space};
\draw [dashed] (0,-2) -- (0,2) node[left]{time};
\draw (-.5,-2) node[left]{$a$} -- (.5,2);
\draw (-2,-.5) -- (2,.5) node[right]{$b$};
\end{tikzpicture}
\end{center}
\caption{A two-dimensional Minkowski diagram containing orthogonal lines $a$ and $b$.} \label{MinOrth}
\end{figure}

Note that, when restricted to only the spatial dimensions, the Minkowskian inner product becomes the three-dimensional Euclidean inner product. Distance and orthogonality behave as usual. When we include the time dimension, we find that the spacetime points which lie a given distance from the origin lie on hyperboloids. When we restrict our diagram to one spatial dimension, orthogonal lines are given a simple geometric description: two lines are orthogonal if they have reciprocal slopes (Figure \ref{MinOrth}). 

A bijection $f:\R^4\to\R^4$ is an \emph{isometry} if it preserves Minkowskian distance; that is, if for all $x,y\in \R^4$, $d_{Minkowski} (x,y)=d_{Minkowski}(f(x),f(y))$. The set of isometries of spacetime, together with the composition opertation, is called the Poincar\'e group.

We may partition the intervals (and lines) of $\R^4$ into three types: \emph{spacelike, timelike,} and \emph{lightlike} intervals. Spacelike intervals can be taken by transformations in the Poincar\'e group to intervals of slope $0$; timelike intervals can be taken to intervals with undefined slope; and lightlike intervals are taken by all transformations to intervals with slope of magnitude $1$. It can be shown that this really does partition $\R^4$.

In drawing a Minkowski diagram, and indeed when describing the movement of objects in the universe, we have been using a fixed choice of coordinate axes. In this coordinate system, there is a (theoretical) object $a$ whose world line goes vertically through the origin. We will call the coordinate system $a$'s \emph{reference frame}, and we will say that $a$ is an \emph{observer}. $a$ is then stationary in its own reference frame. 

The coordinate axes we choose, however, are to some degree arbitrary. We may choose a new set of coordinate axes that we use to describe the movement of objects in spacetime. The world line of some (theoretical) object $b$ runs vertically through the origin of this coordinate system, and we will call the coordinate system the $b$'s reference frame. The change of coordinates from $a$'s reference frame to $b$'s reference frame is called a \emph{reference frame transformation}.

\begin{figure}
\begin{center}
\begin{tikzpicture}[scale=.5]
\draw (-8,0) -- (-2,0) node[right]{$x$};
\draw (-5,-3) -- (-5,3) node[right]{$t$};
\draw [dotted] (-6,-3) -- (-4,3) node[right]{$t'$};
\draw [dotted] (-8,-1) -- (-2,1) node[right]{$x'$};
\draw [->] (-.5,0) -- (.5,0);
\draw (2,0) -- (8,0) node[right]{$x'$};
\draw (5,-3) -- (5,3) node[right]{$t'$};
\draw [dotted] (4,3) node[left]{$t$} -- (6,-3);
\draw [dotted] (8,-1) node[right]{$x$} -- (2,1);
\end{tikzpicture}
\begin{tikzpicture}[scale=.5]
\draw (-8,0) -- (-2,0) node[right]{$x$};
\draw (-5,-3) -- (-5,3) node[right]{$t$};
\draw [dashed] (-8,-3) -- (-2,3);
\draw [dashed] (-8,3) -- (-2,-3);
\draw [->] (-.5,0) -- (.5,0);
\draw (2,0) -- (8,0) node[right]{$x'$};
\draw (5,-3) -- (5,3) node[right]{$t'$};
\draw [dashed] (8,-3) -- (2,3);
\draw [dashed] (8,3) -- (2,-3);
\end{tikzpicture}
\end{center}
\caption{A schematic of a reference frame transformation, in two dimensions. The top diagram indicates the choice of coordinate axes, while the bottom diagram shows that the light cones are invariant under this transformation} \label{transformdiagram}
\end{figure}
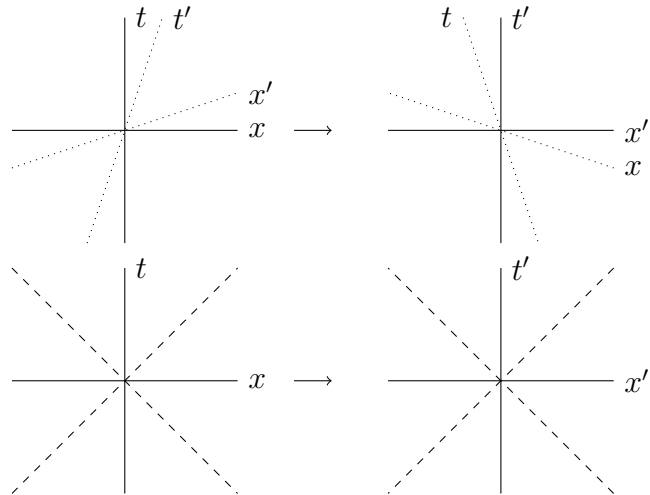

There are some constraints that are imposed on what we may consider a reference frame. We accept the \emph{principle of relativity}, according to which the laws of nature are the same in every reference frame. A consequence of this principle is that reference frame transformations are just those in the Poincar\'e group, and observers are just those objects whose world lines are timelike in every reference frame. What is more, the speed of light is the same in every coordinate system; photons always travel across lightlike intervals. For this reason, light signals are often drawn as a double cone centered at the origin, running along the time axis. Figure \ref{transformdiagram} gives a schematic of a reference frame transformation.

One of the insights of relativity is that some physical properties are relative to which reference frame they are described in. One such property is simultaneity in time. Given two points in spacetime, we say that these points are \emph{simultaneous with respect to an observer $a$} if, in $a$'s reference frame, these points have the same time coordinate. Equivalently, two points are simultaneous with respect to an observer $a$ if the line they define is orthogonal to $a$'s world line. There is no way in which two points can be absolutely simultaneous--simultaneous in all reference frames--because two points which are simultaneous in one reference frame will not be simultaneous in another. 

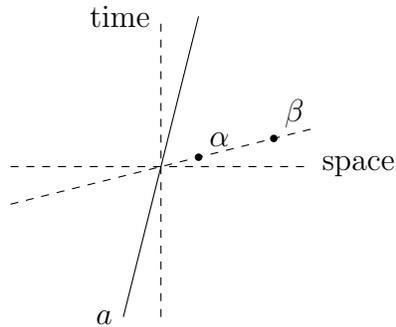
\begin{figure}
\begin{center}
\begin{tikzpicture}[scale=1]
\draw [dashed] (-2,0) -- (2,0) node[right]{space};
\draw [dashed] (0,-2) -- (0,2) node[left]{time};
\draw (-.5,-2) node[left]{$a$} -- (.5,2);
\draw [dashed] (-2,-.5) -- (.5,.125) node[above right]{$\alpha$} -- (1.5,.375) node[above right]{$\beta$} -- (2,.5);
\fill (.5,.125) circle (.05);
\fill (1.5,.375) circle (.05);
\end{tikzpicture}
\end{center}
\caption{A two-dimensional Minkowski diagram containing one observer $a$. The points $\alpha$ and $\beta$ lie on a line orthogonal to $a$, and so they are simultaneous with respect to $a$.} \label{MinSim}
\end{figure}

Not all physical properties depend on which reference frame they are represented in. For instance, if some points $\alpha,\beta$ are separated by a timelike interval in one reference frame, they will be separated by a timelike interval in every reference frame. Because Poincar\'e transformations are isometries, the Minkowskian distance from $\alpha$ to $\beta$ will be the same in every reference frame.

\chapter{Logic Background} \label{LogicChapter}
\section{Introduction}
Usually in mathematics we use a \emph{language} to describe and make arguments about mathematical \emph{structures} that interest us. In logic, the object of study is the \emph{interaction} between these languages, as well as the inference rules we use in deductive reasoning, and the mathematical structures they describe. In this chapter I outline the basic terminology used in mathematical logic, as well as some results that are important, interesting, or especially relevant to this paper.

We start by reviewing the notions of a structure and a language. We next connect the two by using interpretations of terms, and discuss differences between defined and undefined terms. We turn to the basics of model theory: logical implication, elementary equivalence of structures, theories of structures, axiomatizability, and independence of sentences from a theory. Next, we examine the relationship between syntactic and semantic notions of logical consequence, including the classic completeness and soundness theorems. We also look at the decidability of a theory, and methods for proving decidability. Finally, we touch on some important results from model theory: the compactness theorem, and the theorems of Lowenheim and Skolem. Throughout this chapter I use real-closed ordered fields as an example. This is because, first, they are interesting; and second, because they will be of some importance later in this paper.

As was the case in the previous chapter, my description of mathematical logic is not exhaustive and is concentrated on those results which are important to understand the content of this thesis. Any reader interested in a more comprehensive presentation should consult \cite{Manzano}, from which much of this chapter was distilled.

\section{Basics}

\subsection{Structures}

These definitions fall under a field of math called \emph{Universal Algebra}. The motivation is to generalize the way we think about many mathematical objects, whether they be groups, fields, graphs, or ordered sets. Groups are structures, for instance, and a group isomorphism is considered to be an isomorphism in the sense defined below.

\begin{defi} \label{DefiStructure} $\M$ is a \emph{structure} iff $\M=\langle \mathbf{M}, \langle \mathbf{c}_i \rangle_{i \in I} , \langle \mathbf{f}_j \rangle_{j\in J}, \langle \mathbf{R}_k \rangle_{k\in K} \rangle$ such that:

(1) $\mathbf{M}\ne \emptyset$ is the universe of the structure. $I$, $J$, and $K$ are index sets.

(2) $\mathbf{c}_i$ are distinguished elements of $\mathbf{M}$. 

(3) $\mathbf{f}_j:\mathbf{M}^{\mu(j)} \to \mathbf{M}$ are functions, where $\mu: J \to \N \backslash \{0\}$ is a function giving the arity of each $\mathbf{f}_j$.

(4) $\mathbf{R}_k\subset \mathbf{M}^{\delta (k)}$, where $\delta: K \to \N \backslash \{0\}$ is a function giving the arity of each $\mathbf{R}_k$. All structures have equality as a binary relation, even when this is not specified.
\end{defi}

In the remainder of this paper, we denote structures by italics letters like $\M$, $\sN$, and $\mathcal{O}$. Sets, distinguished elements, functions, and relations will be denoted by bold faced letters. This notation will be important later on when we consider the formal languages we use to describe structures. For instance, we may distinguish between a function $\mathbf{f}$ and the function symbol $f$ (in the formal language) that denotes this element.

We define a substructure in the natural way. Here, $\mathbf{f}_j\upharpoonright{(\mathbf{M}^*)^{\mu(j)}}$ denotes the restriction of $\mathbf{f}_j$ to the domain $(\mathbf{M}^*)^{\mu(j)}$.

\begin{defi} $\M^*=\langle \mathbf{M^*}, \langle \mathbf{c}_i^* \rangle_{i \in I} , \langle \mathbf{f}_j^* \rangle_{j\in J}, \langle \mathbf{R}_k^* \rangle_{k\in K} \rangle$ is a substructure of $\M=\langle \mathbf{M}, \langle \mathbf{c}_i \rangle_{i \in I} , \langle \mathbf{f}_j \rangle_{j\in J}, \langle \mathbf{R}_k \rangle_{k\in K} \rangle$ , denoted $\M^* \le \M$, iff: 

(1) $\mathbf{M}^* \subset \mathbf{M}$.

(2) $\mathbf{c}_i^*=\mathbf{c}_i$

(3) $\mathbf{f}_j^*=\mathbf{f}_j\upharpoonright{(\mathbf{M}^*)^{\mu(j)}}$

(4) $\mathbf{R}_k^*=\mathbf{R}_k \cap (\mathbf{M}^*)^{\delta(j)}$
\end{defi}

Now that we have defined structures, we define homomorphisms and isomorphisms between structures. In the following, we let $\M=\langle \mathbf{M}, \langle \mathbf{c}_i \rangle_{i \in I} , $ $\langle \mathbf{f}_j \rangle_{j\in J}, \langle \mathbf{R}_k \rangle_{k\in K} \rangle$ and $\sN=\langle \mathbf{N}, \langle \mathbf{d}_i \rangle_{i \in I} , \langle \mathbf{g}_j \rangle_{j\in J}, \langle \mathbf{S}_k \rangle_{k\in K} \rangle$ be structures.

\begin{defi} The function $\mathbf{h}:\mathbf{M} \to \mathbf{N}$ is a homomorphism from $\M$ to $\sN$ iff, for all $i\in I$, $j \in J$, and $k \in K$:

(1) $\mathbf{h}(\mathbf{c}_i)=\mathbf{d}_i$.

(2) For every $\mathbf{x}_1 ,\dots, \mathbf{x}_{\mu(j)} \in \mathbf{M}$, 

\quad\qquad$\mathbf{h}(\mathbf{f}_i(\mathbf{x}_1,\dots,\mathbf{x}_{\mu(j)}))=
\mathbf{g}_i(\mathbf{h}(\mathbf{x}_1),\dots,\mathbf{h}(\mathbf{x}_{\mu(j)})).$

(3) For every $\mathbf{x}_1 ,\dots, \mathbf{x}_{\delta(k)} \in \mathbf{M}$, if $\langle\mathbf{x}_1,\dots,\mathbf{x}_{\delta(k)}\rangle \in \mathbf{R}_k$, 

\quad\qquad then $\langle\mathbf{h}(\mathbf{x}_1),\dots,\mathbf{h}(\mathbf{x}_{\delta(k)})\rangle \in \mathbf{S}_k$
\end{defi}

When $\mathbf{h}$ is a bijection, and the third condition reads ``\emph{iff}'' instead of ``\emph{if\dots then\dots}'', then $\mathbf{h}$ is an \emph{isomorphism}, and we say $\M$ is isomorphic to $\sN$, denoted $\M \cong \sN$. If there is some structure $\mathcal{O} \le \sN$ such that $\M \cong \mathcal{O}$, then the isomorphism from $\M$ to $\mathcal{O}$ is called an \emph{embedding} of $\M$ into $\sN$.

\subsection{Languages and first-order Logic}

There are many formal languages that are used to describe structures. Here we are only concerned with those which use first-order logic. A language consists of an alphabet of symbols, which can be strung together to make logical sentences. It also contains syntactic rules for when a string of symbols is well formed, that is, when it is meaningful and can be evaluated for its truth value.

All languages share the logical symbols $\neg$, $\land$, $\lor$, $\Rightarrow$, $\Leftrightarrow$, $\forall$, and $\exists$, as well as the non-logical relation symbol $=$. Languages also contain a countably infinite set of variables $x_1,\dots,y_1,\dots,z_1,\dots$ In addition, a language contains the non-logical constant symbols $c_1,\dots,c_i$, function symbols $f_1,\dots,f_j$, and relation symbols $R_1,\dots,R_k$. Each function and relation symbol is assigned an arity. If a language $L$ ``fits'' a structure $\M$, that is, if for each $c_i, f_j, R_k\in L$ there is a corresponding $\mathbf{c}_i, \mathbf{f}_j, \mathbf{R}_k\in \M$, and the arities of the function and relation symbols correspond to the arities of the functions and relations, then we say that $L$ is \emph{adequate} to $\M$. In this paper, we use normal, non-bold letters for the symbols in our formal language.

We distinguish between \emph{terms} and \emph{formulas}. Briefly, a term $\tau$ of a language is a string of symbols that, when interpreted in a structure, will refer to an element of that structure. For instance, in a language $L$ in which $1,2,3,\dots$ are constant symbols and $+$ is a 2-ary function, then $4$, $5+2$, and $7+x$ are all terms. A formula $\phi$ of a language is a string of symbols that, when interpreted in a structure, will have a truth value. In the same language $L$, $3+4=7$, $\forall x(x=2)$, and $\exists y(y+x=3)$ are all formulas.

Using induction, one can rigorously define the rules by which a string of symbols is well-formed. I skip this here, and instead give some examples. The strings $f(x,y,z)=c$, $\forall x(x=y) \Rightarrow \forall x (x+x=y)$, and $\exists x,y (R(x,x) \land \neg R(x,y)$ are well-formed, provided that the arity of $f$ is three, that the language contains the 2-ary function symbol $+$, and $R$ is a 2-ary relation.

The strings $\neg \lor x$, $\forall f( f(x)=x)$, and $\exists c (c\land x=y)$ are not well formed. The second and are problematic because in first-order logic, the only symbols that can be quantified over  are variables. The third is also problematic because, when interpreted in a structure, $c$ is a term, not a formula, and it will not have a truth value when interpreted in a structure. Hence, $c\land x=y$ is meaningless.

We call variables \emph{free} if they are not quantified over. If a variable is quantified over, it is \emph{bound}. For instance, in the string $f(x,y,z)=c$, the variables $x$, $y$, and $z$ are free. In $\forall x(x=y) \Rightarrow \forall x (x+x=y)$, the variable $y$ is free, while $x$ is bound. When $x_1,\dots,x_n$ are among the free variables of a formula $\phi$, we write $\phi(x_1,\dots,x_n)$. When a formula $\phi$ contains no free variables, we call $\phi$ a \emph{sentence}.

\subsection{Connecting a language with a structure}

When we interpret a language in a structure, we are essentially assigning a meaning to each of the symbols in our language. Some symbols, like $\neg$, $\forall$, and $=$ have a meaning that does not vary depending on interpretation. Our constant, function, and relation symbols will stay fixed in meaning for all interpretations in a single structure. Our variable symbols can vary in meaning across each interpretation, even if these interpretations are in the same structure. More formally now:

\begin{defi} An \emph{assignment} of a language $L$'s variables to a structure $\M$ is a function $\mathcal{I}:\mathcal{V}\to \mathbf{M}$, where $\mathcal{V}$ is the set of variables in $L$. For an assignment that assigns the variable $x$ to the element $\mathbf{x}$, we write $\mathcal{I}_x^\mathbf{x}$.
\end{defi}

We now inductively define the interpretation of formulas in a structure. Intuitively, a structure $\M$ \emph{satisfies} a formula $\phi$ when, given an assignment of the free variables in $\phi$, $\phi$ holds true about $\M$.

\begin{defi} Given a structure $\M$ and an assignment $\mathcal{I}$, we define the interpretation $\M\mathcal{I}$ as follows:

$\M\mathcal{I}(x)=\mathcal{I}(x)$,

$\M\mathcal{I}(f_j \tau_1,\dots,\tau_{\mu(j)})=\mathbf{f}_j(\M\mathcal{I}(\tau_1),\dots,\M\mathcal{I}(\tau_{\mu(j)})$,

$\M\mathcal{I}~\rm{sat}~R_k\tau_1\dots\tau_{\delta(k)}~iff~\langle\M\mathcal{I}(\tau_1),\dots,\M\mathcal{I}(\tau_{\delta(k)})\rangle\in\mathbf{R}_k$,

$\M\mathcal{I}~\rm{sat}~\tau_1=\tau_2~iff~\M\mathcal{I}(\tau_1)=\M\mathcal{I}(\tau_2)$,

$\M\mathcal{I}~\rm{sat}~\neg\phi~iff~\rm{not}~\M\mathcal{I}~\rm{sat}~\phi$,

$\M\mathcal{I}~\rm{sat}~(\phi_1\land\phi_2~\rm{/}~\phi_1\lor\phi_2~\rm{/}~\phi_1\Rightarrow\phi_2~\rm{/}~\phi_1\Leftrightarrow\phi_2)~iff$

\qquad\qquad$\M\mathcal{I}~\rm{sat}~\phi_1~(\rm{and/or/implies/iff})~\M\mathcal{I}~\rm{sat}~\phi_2$,

$\M\mathcal{I}~\rm{sat}~\forall x\phi~iff~\rm{for~all}~\mathbf{x}\in\mathbf{M}: \M\mathcal{I}_x^\mathbf{x}~\rm{sat}~\phi$,

$\M\mathcal{I}~\rm{sat}~\exists x\phi~iff~\rm{for~some}~\mathbf{x}\in\mathbf{M}: \M\mathcal{I}_x^\mathbf{x}~\rm{sat}~\phi$.
\end{defi}

Note that whether $\M\mathcal{I}$ satisfies a sentence $\phi$ does not depend on any specific choice of $\mathcal{I}$. If $\phi$ is a formula, we say that $\M$ \emph{models} $\phi$ iff, for any assignment $\mathcal{I}$, $\M\mathcal{I}~\rm{sat}~\phi$. Thus, if $\phi$ is a sentence, $\M$ models $\phi$, denoted $\M\models\phi$, iff there exists an assignment $\mathcal{I}$ such that $\M\mathcal{I}~\rm{sat}~\phi$. If $\Gamma$ is a set of sentences such that $\M\models\phi$ for each $\phi\in\Gamma$, we say that $\M$ models $\Gamma$, again denoted $\M\models\Gamma$. We let $\rm{Mod}(\Gamma)=\{\M | \M\models\Gamma\}$.

Some simple examples: Any structure $\M$ models the formula $x=x$. Any structure $\M$ has $\M\models\forall x (x=x)$. Any group $\mathcal{G}$ has $\mathcal{G}\models\exists x \forall y (x*y=y)$. The ordered reals $\langle\R,\mathbf{<}\rangle$ are such that $\langle\R,\mathbf{<}\rangle\models\forall x y (x<y \Rightarrow \exists z(x<z<y))$.

\subsection{Definability}

In a language $L$, the constant, function, and relation symbols we begin with are undefined; they acquire meaning only under an interpretation in a structure. We can also define new relation symbols using the undefined symbols. These new symbols will acquire meaning as the undefined symbols they are built out of acquire meaning. A term defined within a language is then, to some extent, just short hand for a longer string of symbols. We now treat this formally.

\begin{defi}
Let $\M=\langle \mathbf{M}, \langle \mathbf{c}_i \rangle_{i \in I} , \langle \mathbf{f}_j \rangle_{j\in J}, \langle \mathbf{R}_k \rangle_{k\in K} \rangle$ be a structure. The relation $\mathbf{R} \subset \mathbf{M}^n$ is \emph{definable} iff there exists a formula $\phi\langle x_1\dots x_n\rangle$ such that $\mathbf{R}=\{\langle\mathbf{x}_1,\dots,\mathbf{x}_n\rangle\in\mathbf{M}^n | \M\mathcal{I}_{x_1\dots x_n}^{\mathbf{x}_1\dots \mathbf{x}_n}~\rm{sat}~\phi\}$.
\end{defi}
Consider the following examples: In the structure $\langle\R,\mathbf{<}\rangle$, letting $\phi(x,y)\equiv x=y\lor x<y$ shows that the binary relation $\le$ is definable. This is because, in $\langle\R,\mathbf{<}\rangle$, $\mathbf{\le}=\{\langle x,y \rangle \in \R^2 | \mathbf{x=y}~\rm{or}~\mathbf{x<y}\}$.

In the structure $\mathcal{C}_{\R^3}=\langle \mathbf{\R}^3 \cup \mathbf{\R}, \mathbf{0}, \mathbf{1}, \mathbf{|\cdot |}, \mathbf{+}, \mathbf{*}, \mathbf{<}\rangle$, where $\langle\mathbf{\R},\mathbf{+},\mathbf{*},\mathbf{<}\rangle$ is the usual ordered field and $\mathbf{|\cdot |}:\mathbf{\R}^3 \times \mathbf{\R}^3 \times\mathbf{\R}^3 \to \mathbf{\R}$ is the determinant of the matrix formed by three vectors in $\mathbf{\R}^3$, we can define the 3-ary ``collinear'' relation $\mathbf{\rm{Col}}$ by means of the formula 
\begin{equation} \phi(x,y,z) \equiv 
\begin{vmatrix} x_1 & x_2 & x_3 \\
y_1 & y_2 & y_3 \\
z_1 & z_2 & z_3 \end{vmatrix} = 0.
\end{equation}
This works because $\mathbf{x},$ $\mathbf{y},$ and $\mathbf{z}$ are collinear iff $\mathcal{C}_{\R^3}\mathcal{I}_{x,y,z}^\mathbf{x,y,z}~\rm{sat}~\phi(x,y,z)$; that is, in $\mathcal{C}_{\R^3}$, $\mathbf{\rm{Col}} =\{ \mathbf{x,y,z} | \text{the determinant of } (\mathbf{x,y,z})=0\}$.

We note that, in this case, $\mathbf{|\cdot |, +, *, <}$ are functions and relations in our structure, but $|\cdot |, +, *, <$ are also the symbols in our language that denote these functions and relations under an interpretation. How these symbols, and symbols like these, are being used in this paper will be apparent from context.

\subsection{Many-Sorted Structures and Languages}

Sometimes it is convenient to give a structure with more than one universe. A structure with $n$ pairwise disjoint universes is said to be \emph{$n$-sorted} and is written: $\mathbf{\langle \langle X_n,\rangle_{n\in N},\langle \mathbf{c}_i \rangle_{i \in I} , \langle \mathbf{f}_j \rangle_{j\in J}, \langle \mathbf{R}_k \rangle_{k\in K} \rangle}$. For each function $\mathbf{f_j}$ the entries of each argument of $\mathbf{f_j}$ may draw from only one $\mathbf{X_n}$. The same is true for each relation $\mathbf{R_k}$.

A language adequate to a $n$-sorted structure will be an $n$-sorted language in which there are $n$ sorts of variables, which each denote elements from a single universe. Typically, these different sorts of variables will be written with different kinds of letters; for instance, Roman for one sort, and Greek for another. In giving a language, it is specified which sort of variables may serve as entries of each argument of the function and relation symbols in that language.

It is a simple task to reduce a $n$-sorted structure to a single-sorted one, by taking the new universe of the structure to be the union of the $n$ sorted universes and adding in $n$ 1-ary predicates corresponding to the $n$ sorted universes. A similar strategy reduces a $n$-sorted language to a single-sorted one. Expressing structures and languages as many-sorted is then only a matter of notational simplicity. Two $n$-sorted structures $\mathbf{\langle \langle X_n,\rangle_{n\in N},\langle \mathbf{c}_i \rangle_{i \in I} , }$ $\mathbf{\langle \mathbf{f}_j \rangle_{j\in J}, \langle \mathbf{R}_k \rangle_{k\in K} \rangle}$ and $\mathbf{\langle \langle X'_n,\rangle_{n\in N},\langle \mathbf{c'}_i \rangle_{i \in I} , \langle \mathbf{f'}_j \rangle_{j\in J}, \langle \mathbf{R'}_k \rangle_{k\in K} \rangle}$ are isomorphic iff the corresponding single-sorted structures are isomorphic. Equivalently, they are isomorphic if there exist bijections $\phi_n:\mathbf{X_n\to X'_n}$ which appropriately preserve constants, functions, and relations.

\section{Model Theory}
\subsection{Logical Implication and Elementary Equivalence}

Using the notion of a model of a set of sentences, we can define some important basic logical notions.

\begin{defi}
We say that a formula $\phi$ \emph{logically implies} a formula $\psi$, denoted $\phi\models\psi$, iff any structure $\M$ that models $\phi$ also models $\psi$.
\end{defi}

Note that the symbol $\models$ is used both when a structure models a sentence, and also when one formula logically implies another. When every structure models a formula $\phi$, we say that $\phi$ is \emph{logically valid}, denoted $\models\phi$. For instance, $\models x=x$. When $\phi\models\psi$ and $\psi\models\phi$, we say that $\phi$ and $\psi$ are \emph{logically equivalent.} This is equivalent to writing $\models\phi\Leftrightarrow\psi$.

\begin{defi}
Given a structure $\M$, the \emph{theory of $\M$} $\rm{Th}(\M)$ is the set of all sentences $\phi$ such that $\M\models\phi$. Given a class of structures $\Sigma$, the theory of $\Sigma$ is the set of all sentences $\phi$ such that, for each $\M\in\Sigma$, $\M\models\phi$.
\end{defi}

Thus, if $\mathcal{R}=\langle\R,0,1,+,*,<\rangle$ is the real ordered field, then $\rm{Th}(\mathcal{R})$ is the set of sentences that hold true about $\mathcal{R}$. The theory of groups is the set of sentences that hold true about all structures that are groups.

\begin{defi}
Given two structures $\M$ and $\sN$ and a language $L$ that is adequate to them both, we say that $\M$ is \emph{elementary equivalent} to $\sN$, denoted $\M\equiv\sN$, if for every sentence $\phi$ written in $L$, $\M\models\phi$ iff $\sN\models\phi$.
\end{defi}

Two structures are elementary equivalent when they are indistinguishable from the perspective of first-order logic. That is, any thing we can say in first-order logic about one also holds for the other. We now consider an important example.

We consider $\rm{Mod}(\rm{Th}(\mathcal{R}))$, the set of structures which satisfy all the sentences true of $\mathcal{R}$. Then $\M\in\rm{Mod}(\rm{Th}(\mathcal{R}))$ iff $\M\equiv\mathcal{R}$. We call these structures the \emph{real-closed ordered fields}. Obviously, $\mathcal{R}$ is a real-closed ordered field. We'll see later that there is a real-closed ordered field for each infinite cardinality, but for now let us consider one specific real-closed ordered field. The \emph{algebraic numbers} $\mathbf{A}$ are those numbers that are solutions to one-variable polynomials $a_0+a_1x+a_2x^2+\dots+a_ix^i$, where each $a_j$, $1\le j\le i$, is rational. The ordered field $\langle \mathbf{A},\mathbf{0},\mathbf{1},+,*,<\rangle$ of $\mathbf{A}$ is real-closed. It is also countable.

This highlights the degree to which first-order logic is not able to capture all the properties of a structure. In order to find a set of sentences $\Gamma$ such that, for each $\M\models\Gamma$, $\M\cong\mathcal{R}$, we would have to move into \emph{second-order} logic, in which quantifiers can range over sets as well as elements. In this more complex logic, we could formulate a sentence about the Dedekind-completeness of the reals. If a theory $T$ is such that, for all structures $\M$ and $\sN$, if $\M\models T$ and $\sN\models T$, then $\M\cong\sN$, then we say that $T$ is a \emph{categorical} theory. We just noted, then, that $\rm{Th}(\mathcal{R})$ is not categorical.

Another property that is important to theories is that of \emph{completeness}. A theory $T$ in a language $L$ is complete iff, for every sentence $\phi$ written in $L$, either $\phi\in T$ or $\neg\phi\in T$. The theory of groups is not complete: some groups are abelian and some groups are not abelian. $\rm{Th}(\mathcal{R})$ is complete because, for each sentence $\phi$ in the language of ordered fields, $\mathcal{R}$ satisfies $\phi$ or satisfies $\neg\phi$.

\subsection{Axiomatizability and Independence}

Given a theory $T$ in a language $L$, we may ask if there is some decidable set of sentences $\Gamma$ such that for any structure $\M$, $\M\models\Gamma$ iff $\M\models T$. If there is, we say that $T$ is \emph{axiomatizable}. If there is some $\Gamma$ that is finite that axiomatizes $T$, then $T$ is \emph{finitely axiomatizable}.

For a set to be decidable, there must be some algorithm that tells us, for each sentence $\phi$ written in $L$, if it is true that $\phi\in\Gamma$. An algorithm is a finite set of well-defined rules that given an appropriate question (in this case ``is $\phi$ in $\Gamma$?''), will always respond after a finite number of steps with an answer, and this answer will always be correct.

For instance: consider $\rm{Th}(\mathcal{G})$, the theory of groups, written in a first-order language $L$ adequate to $\mathcal{G}$. This is finitely axiomatizable by three sentences, those which state in $L$ that (1) multiplication is associative, (2) an identity element exists, and (3) inverses exist.

It is simple to give a finite axiomatization for the class of ordered fields. If we add to these axioms a sentence $\forall x(x>0 \Rightarrow \exists y (y^2 =x))$, we arrive at a (finite) axiomatization of  a class of structures called Euclidean fields. If we add to these axioms a set of sentences $\{ \phi_i | i\in\N\}$, where $\phi_i \equiv \forall a_{2i+1},a_{2i},\dots,a_1,a_0\exists x (a_{2i+1}x^{2i+1}+a_{2i}x^{2i}+\dots+a_1x+a_0 = 0$ is a statement asserting that all $2i+1$-degree polynomials have a root in the field, then we arrive at an (infinite) axiomatization of the real-closed ordered fields (see \cite{Tarski}).

One important property we can show about an axiomatization is that of \emph{independence}. We say that a sentence $\phi$ is independent of a set of sentences $\Gamma$ iff $\Gamma$ does not logically imply $\phi$, written $\Gamma\not\models\phi$. If a sentence in an axiomatization is not independent of the others, it could be eliminated in favor of a simpler axiomatization without changing which structures are being axiomatized.

Completeness becomes especially interesting when we consider a theory that is given by only an axiomatization. It is often nontrivial whether such an axiomatization generates a complete theory.

\section{Other methods in logic}

\subsection{Syntactic Methods}

Thus far we have been concerned with how interpretation in a model can help us make sense of some basic logical notions, like implication and independence. This is a \emph{semantic} approach, since the properties of formulas in our language depend on their meanings when interpreted in structures. In this subsection we explore a second approach, which is \emph{syntactic} rather than semantic. We define a deductive calculus, which is a finite set of rules by which we can deduce a formula $\phi$ from some set of formulas $\Gamma$ we take as a hypothesis. The deductive steps we make depend not on the meaning of the formulas in $\Gamma$, but instead on only their syntactic properties. When we can deduce $\phi$ from $\Gamma$, we write $\Gamma\vdash\phi$. If $\Gamma=\emptyset$, we write $\vdash\phi$.

I will not present an entire deductive calculus, but instead give some of its rules as examples. I follow \cite{Manzano} for the rules of the calculus.

Rule HI: If $\phi \in \Gamma$, then $\Gamma \vdash \phi$.

Rule NC: If $\Gamma \cup \{\phi\} \vdash \psi$ and $\Gamma \cup \{\phi\} \vdash \neg\psi$, then $\Gamma\vdash\neg\phi$.

Rule IDC: If $\Gamma\vdash\phi$, then $\Gamma\vdash\phi\lor\psi$.

The deductive calculus is supposed to capture and formalize all the logical steps we make in writing proofs, or at least in those proofs that don't require machinery unavailable to first-order logic, like quantifying over sets. In a sense, we may consider the model-theoretic notion of logical implication ($\models$) as giving the ``true'' relationship between sentences written in our proof language, because it is grounded by interpretations of these sentences in the mathematical objects we study. It then becomes important to know about the relationship between $\models$ and $\vdash$. If we prove something, then do we know that it is ``actually'' the case? And are there ``facts'' about mathematical objects that are not provable? There are two important theorems that speak to this.

\begin{theo}{Soundness Theorem:}
For every set of first-order formulas $\Phi$ and every first-order formula $\psi$, if $\Phi\vdash\psi$, then $\Phi\models\psi$.
\end{theo}

\begin{theo}{Completeness Theorem:} \label{ComTheo}
For every set of first-order formulas $\Phi$ and every first-order formula $\psi$, if $\Phi\models\psi$, then $\Phi\vdash\psi$.
\end{theo}

The soundness theorem tells us that, if we can deduce $\psi$ from $\Phi$, then $\Phi$ logically implies $\psi$. In a sense: our rules of deduction will not lie to us. The completeness theorem tells us that, if $\Phi$ logically implies $\psi$, then we can use our rules of deduction to deduce $\psi$ from $\Phi$. In a sense: there are no logically valid statements that our rules of deduction are blind to. So, when using first-order logic, we can prove any (model-theoretically) true statement, and any statement we prove will be (model-theoretically) true.

If a set of sentences $\Gamma$ is such that there is no formula $\phi$ such that $\Gamma\vdash\phi$ and $\Gamma\vdash\neg\phi$, then $\Gamma$ is \emph{logically consistent}. If $\Gamma$ has a model, then $\Gamma$ is \emph{satisfiable}. The soundness and completeness theorems tell us that, in first-order logic, the two conditions are equivalent.

``Completeness'' here should not be confused with ``completeness'' elsewhere in this paper. Here, ``completeness'' refers to a property of a logical system. Elsewhere, it refers to the property of a theory.

\subsection{Decidability}

One potential advantage to investigating a theory $T$ in the terms of first-order logic is that $T$ may be decidable: there may be an algorithm that tells us, for each sentence $\phi$ written in $L$, if $\phi\in T$. If $T$ is decidable, we could program a computer to tell us what sentences are contained within $T$.

While decidability is not the same property as completeness, they are connected:

\begin{prop} 
\label{DecidabilityCriterion} If a theory is axiomatizable and complete, then it is decidable.
\end{prop}
Also, if a theory is decidable, it is axiomatizable. In this case, the theory is an axiomatization of itself.

One common method for proving the decidability of a theory is \emph{quantifier elimination}. A theory $T$ admits quantifier elimination iff, for every sentence $\phi\in T$, there exists a sentence $\phi*$ such that $\phi*$ has no quantifiers, and $T\vdash \phi \Leftrightarrow \phi*$. If we have an algorithm for translating sentences in $T$ into quantifier-free versions, and if we can show that the set of quantifier-free sentences of $T$ is decidable, then we have shown that $T$ as a whole is decidable.

A notable application of the method of quantifier elimination is in Tarski's proof in \cite{Tarski2} of the decidability of real-closed ordered fields. In section \ref{secEucGeo} we look at a second example of a decidable theory: the first-order theory of Euclidean geometry.

\section{Interesting results in model theory}

In this section, I state without proof some interesting results from model theory, as well as some of their consequences. They are not strictly relevant to the main argument of this paper.

\begin{theo}{Compactness:} A set of sentences $\Gamma$ has a model iff every finite subset of $\Gamma$ has a model.
\end{theo}

Using the compactness theorem, we can construct \emph{nonstandard models} of the theories of well-known structures, like $\langle\mathcal{R},\mathbf{-},\mathbf{||}\rangle$, the real ordered field along with the negation and absolute value functions. In one such model of $\rm{Th}(\langle\mathcal{R},\mathbf{-},\mathbf{||}\rangle)$, called the hyperreals, there exist infinite and infinitesimal numbers. Even so, this structure is elementary equivalent to $\langle\mathcal{R},\mathbf{-},\mathbf{||}\rangle$. If we take a reduct of this structure--throw out of the structure the functions denoted by the function symbols $-$ and $||$--then we arrive at a nonstandard real closed ordered field.

\begin{theo}{L\"owenheim-Skolem:} Let $\Gamma$ be a set of sentences of cardinality $\gamma$. If $\Gamma$ has a model of infinite cardinality $\alpha>\gamma$, then for every $\beta\ge\omega$ such that $\omega\le\beta\le\alpha$, $\Gamma$ has a model of cardinality $\beta$. 
\end{theo}

\begin{theo}{Upward L\"owenheim-Skolem:} Let $\Gamma$ be a set of sentences of cardinality $\gamma$. If $\Gamma$ has a model of cardinality $\alpha\ge\omega$, then for every $\beta\ge\omega,\gamma$, $\Gamma$ has a model of cardinality $\beta$.
\end{theo}

These two theorems emphasize the extent to which first-order theories can fail to be categorical. For instance, we saw earlier that $\rm{Th}(\mathcal{R})$ can be axiomatized by a countably infinite set of sentences $\Gamma$. Since $\Gamma$ has a model of cardinality $\mathfrak{c}>\omega=|\rm{Th}(\mathcal{R})|$, from the L\"owenheim-Skolem theorems it follows that $\rm{Th}(\mathcal{R})$ has a model of \emph{every infinite cardinality}! So the algebraic numbers and the hyperreals are just two of many, many real closed fields.

\chapter{A simple axiomatization of special relativity, with only slower-than-light observers} \label{STLChapter}
\section{Introduction}
In this chapter, we develop an axiom system \textsc{SimpleRel} for special relativity. There are a number of axiomatizations of special relativity in the literature, including a well-developed one put forth by Andr\'eka et. al. in \cite{Renyi}. The motivations behind \textsc{SimpleRel} are unique. It is written in a language with only two undefined symbols. Unlike the axiom system in the work of Andr\'eka et. al., it does not include a field of quantities to describe the physical objects in its universe. It ontologically commits to the existence of only physical objects, and only two kinds of objects: observers and signals. Finally, its axioms were designed to be empirically testable and confirmed by a strong body of evidence.

\textsc{SimpleRel} is a significant revision of an axiomatization made by Ax in \cite{Ax}. Ax made the assumption that all observers travel slower than light. We find that this statement is not independent of the other axioms in \textsc{SimpleRel}, but instead follows from assumptions on the isotropy of space. We prove this using model-theoretic methods: we characterize all the models of \textsc{SimpleRel}, and show that this statement holds true in all these models. By the Completeness theorem \ref{ComTheo}, this statement is provable from \textsc{SimpleRel}. Thus, by using model theory, we are able to show that this statement is provable from \textsc{SimpleRel} in first-order logic \emph{without actually giving a first-order proof for it.}

In the course of his axiomatization, Ax makes use of a finite axiomatization of Euclidean geometry in three dimensions. We use an infinite axiomatization instead. Because of this change, we have that the models of \textsc{SimpleRel} are all isomorphic to structures defined over real-closed ordered fields. Ax can only show that all the models of his axiom system are isomorphic to structures defined over Euclidean fields, a wider class of fields. A consequence of this change is that \textsc{SimpleRel} is complete and decidable. 

It is possible to replace the axiom system $\mathfrak{E}$ in section \ref{secEucGeo} with the finite system of axioms that Ax uses. In this case, our representation theorem in section \ref{RepThmSTLSection} still holds, except that the models of this revised version of \textsc{SimpleRel} are all isomorphic to structures defined over Euclidean fields.

\section{An Axiomatization of Euclidean Geometry}
\label{secEucGeo}

We begin by giving a first-order axiom system for Euclidean Geometry in three dimensions, from \cite{Tarski}. Call this axiom system $\mathfrak{E}$. $\mathfrak{E}$ is written in a language with two non-logical predicates: the 3-ary relation $\be$ and the 4-ary relation $\ed$. $\be(xyz)$ is interpreted: ``$y$ is between $x$ and $z$.'' $\ed(xy,zw)$ is interpreted: ``The distance between $x$ and $y$ is the same as the distance between $z$ and $w$.''

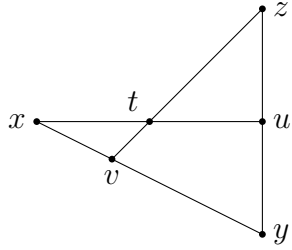
\begin{figure}
\begin{center}
\begin{tikzpicture}[scale=1.5]
\draw (2,0) node[right] {$u$} -- (1,0) node[above left] {$t$} -- (0,0) node[left] {$x$} -- (.667,-.333) node[below] {$v$} -- (2,-1) node[right] {$y$} -- (2,1) node[right] {$z$} -- (.667,-.333);
\fill (2,0) circle (0.03);
\fill (1,0) circle (0.03);
\fill (0,0) circle (0.03);
\fill (.667,-.333) circle (0.03);
\fill (2,-1) circle (0.03);
\fill (2,1) circle (0.03);
\end{tikzpicture}
\caption{\textsc{Axiom 7 (Pasch's axiom)}} \label{PaschAx}
\end{center}
\end{figure}

\medbreak
\noindent \textsc{Axiom 1 (Identity axiom for betweenness)}: 
\begin{equation}
\forall x y (\be(xyx) \Rightarrow (x=y)).
\end{equation}
\textsc{Axiom 2 (Transitivity axiom for betweenness)}: 
\begin{equation}
\forall xyzu(\be(xyu)\land\be(yzu)\Rightarrow\be(xyz)).
\end{equation}
\textsc{Axiom 3 (Connectivity axiom for betweenness)}: 
\begin{equation}
\forall xyzu(\be(xyz) \land\be(xyu)\land(x\ne y) \Rightarrow \be(xzu)\lor\be(xuz))
\end{equation}
\textsc{Axiom 4 (Reflexivity axiom for equidistance)}:
\begin{equation}
\forall xy(\ed(xy,yx))
\end{equation}
\textsc{Axiom 5 (Identity axiom for equidistance)}:
\begin{equation}
\forall xyz(\ed(xy,zz)\Rightarrow(x=y))
\end{equation}
\textsc{Axiom 6 (Transitivity axiom for equidistance)}:
\begin{equation}
\forall xyzuvw (\ed(xy,zu) \land \ed(xy,vw) \Rightarrow \ed(zu,vw))
\end{equation}
\textsc{Axiom 7 (Pasch's axiom)} (Figure \ref{PaschAx}):
\begin{equation}
\forall txyzu [\be(xtu) \land \be(yuz) \Rightarrow \exists v(\be(xvy) \land \be(ztv))]
\end{equation}

\begin{figure}
\begin{center}
\begin{tikzpicture}[scale=1]
\draw (2,2) node[right]{$x$} -- (1,1) node[left]{$y$} -- (0,0) node[left]{$w$} -- (2,0) node[below]{$t$} -- (4,0) node[right]{$v$} -- (3,1) node[right]{$z$} -- (2,2);
\draw (2,2) -- (2,1) node[above left]{$u$} -- (2,0);
\draw (1,1) -- (3,1);
\fill (0,0) circle (0.05);
\fill (2,0) circle (0.05);
\fill (4,0) circle (0.05);
\fill (1,1) circle (0.05);
\fill (2,1) circle (0.05);
\fill (3,1) circle (0.05);
\fill (2,2) circle (0.05);
\end{tikzpicture}
\caption{\textsc{Axiom 8 (Euclid's axiom)}} \label{EuclidAx}
\end{center}
\end{figure}
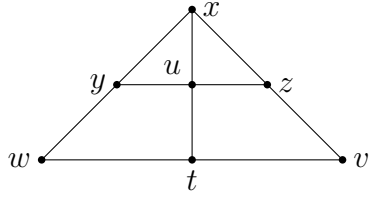

\noindent \textsc{Axiom 8 (Euclid's axiom)} (Figure \ref{EuclidAx}):
\begin{equation}
\begin{split}
\forall txyzu [\be(xut) \land \be(yuz) \land (x \ne u) \\ 
\Rightarrow \exists vw (\be(xzv) \land \be(xyw) \land \be(vtw))]
\end{split}
\end{equation}

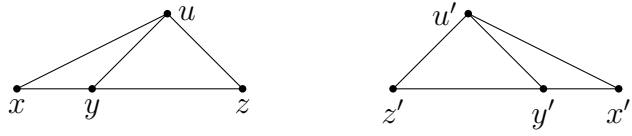
\begin{figure}
\begin{center}
\begin{tikzpicture}[scale=1]
\draw (0,0) node[below]{$x$} -- (1,0) node[below]{$y$} -- (3,0) node[below]{$z$} -- (2,1) node[right]{$u$} -- (0,0);
\draw (1,0) -- (2,1);
\draw (5,0) node[below]{$z'$} -- (7,0) node[below]{$y'$} -- (8,0) node[below]{$x'$} -- (6,1) node[left]{$u'$} -- (5,0);
\draw (6,1) -- (7,0);
\fill (0,0) circle (0.05);
\fill (1,0) circle (0.05);
\fill (3,0) circle (0.05);
\fill (5,0) circle (0.05);
\fill (7,0) circle (0.05);
\fill (8,0) circle (0.05);
\fill (2,1) circle (0.05);
\fill (6,1) circle (0.05);
\end{tikzpicture}
\caption{\textsc{Axiom 9 (Five-segment axiom)}} \label{FSAx}
\end{center}
\end{figure}

\noindent\textsc{Axiom 9 (Five-segment axiom)} (Figure \ref{FSAx}):
\begin{equation}
\begin{split}
\forall xx'yy'zz'uu'(\ed(xy,x'y')\land \ed(yz,y'z') \land \ed(xu,x'u') \land \ed(yu,y'u') \\ \land \be(xyz) \land \be(x'y'z') \land (x\ne y) \Rightarrow \ed(zu,z'u'))
\end{split}
\end{equation}

\begin{figure}
\begin{center}
\begin{tikzpicture}[scale=1]
\draw (0,0) node[above]{$x$} -- (1,0) node[above]{$y$} -- (3.236,0) node[above]{$z$};
\draw (5,-1) node[left]{$u$} -- (6,1) node[right]{$v$};
\fill (0,0) circle (0.05);
\fill (1,0) circle (0.05);
\fill (3.236,0) circle (0.05);
\fill (5,-1) circle (0.05);
\fill (6,1) circle (0.05);
\end{tikzpicture}
\caption{\textsc{Axiom 10 (Axiom of segment construction)}} \label{SegAx}
\end{center}
\end{figure}
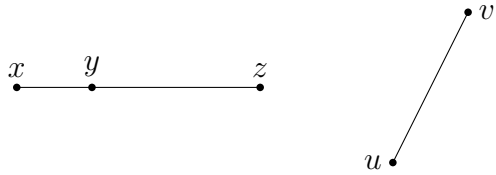

\noindent\textsc{Axiom 10 (Axiom of segment construction)} (Figure \ref{SegAx}):
\begin{equation}
\forall xyuv \exists z(\be(xyz) \land \ed(yz,uv))
\end{equation}

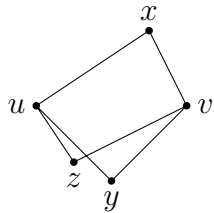
\begin{figure}
\begin{center}
\begin{tikzpicture}[scale=1]
\draw (0,0) node[left]{$u$} -- (1.5,1) node[above]{$x$} -- (2,0) node[right]{$v$} -- (1,-1) node[below]{$y$} -- (0,0) -- (.5,-.75) node[below]{$z$} -- (2,0);
\fill (0,0) circle (0.05);
\fill (2,0) circle (0.05);
\fill (1.5,1) circle (0.05);
\fill (.5,-.75) circle (0.05);
\fill (1,-1) circle (0.05);
\end{tikzpicture}
\caption{\textsc{Axiom 11 (Dimension at least three)} (drawn in three dimensions)} \label{MinDimAx}
\end{center}
\end{figure}

\noindent\textsc{Axiom 11 (Dimension at least three)} (Figure \ref{MinDimAx}):
\begin{equation}
\begin{split}
\exists xyzuv(\ed(xu,xv) \land \ed(yu,yv) \land \ed(zu,zv) \land (u\ne v) \\ \land \neg\be(xyz) \land \neg \be(yzx) \land \neg\be(zxy))
\end{split}
\end{equation}

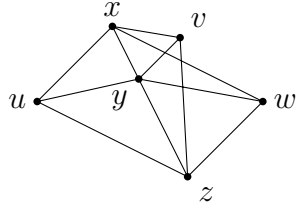
\begin{figure}
\begin{center}
\begin{tikzpicture}[scale=1]
\draw (0,0) node[left]{$u$} -- (1,1) node[above]{$x$} -- (3,0) node[right]{$w$} -- (2,-1) node[below right]{$z$} -- (0,0) -- (1.35,.3) node[below left]{$y$} -- (1,1) -- (1.9,.85) node[above right]{$v$} -- (2,-1) -- (1.35,.3) -- (3,0);
\draw (1.35,.3) -- (1.9,.85);
\fill (0,0) circle (0.05);
\fill (1,1) circle (0.05);
\fill (3,0) circle (0.05);
\fill (2,-1) circle (0.05);
\fill (1.35,.3) circle (0.05);
\fill (1.9,.85) circle (0.05);
\end{tikzpicture}
\caption{\textsc{Axiom 12 (Dimension at most three)} (drawn in three dimensions)} \label{MaxDimAx}
\end{center}
\end{figure}

\noindent\textsc{Axiom 12 (Dimension at most three)} (Figure \ref{MaxDimAx}):
\begin{equation}
\begin{split}
\forall xyzuvw((\ed(xu,xv)\land\ed(xu,xw)\land\ed(yu,yv)\land\ed(yu,yw) \\
\land\ed(zu,zv)\land\ed(zu,zw)\land\neg\be(uvw)\land\neg\be(vuw) \\
\land\neg\be(uwv))\Rightarrow(\be(xyz)\lor\be(yxz)\lor\be(xzy)))
\end{split}
\end{equation}

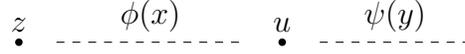
\begin{figure}
\begin{center}
\begin{tikzpicture}[scale=1]
\draw [dashed] (.5,0) -- node[above]{$\phi(x)$} (3,0);
\draw [dashed] (4,0) -- node[above]{$\psi(y)$} (6,0);
\fill (3.5,0) circle (0.05) node[above]{$u$};
\fill (0,0) circle (0.05) node[above]{$z$};
\end{tikzpicture}
\caption{\textsc{Axiom schema 13 (Elementary continuity axioms)}} \label{ContAx}
\end{center}
\end{figure}

\noindent\textsc{Axiom schema 13 (Elementary continuity axioms)} (Figure \ref{ContAx}): 
All sentences of the form
\begin{equation}
\begin{split}
\forall w_1\dots w_n (\exists z\forall xy(\phi(x)\land\psi(y)\Rightarrow\be(zxy)) \\ 
\Rightarrow\exists u\forall xy(\phi(x)\land\psi(y)\Rightarrow\be(xuy))),
\end{split}
\end{equation}
where $\phi$ is a formula such that the variables $x,w_1,\dots,w_n$ may occur free, but $y,z,u$ do not occur free, and where $\psi$ is a formula such that the variables $y,w_1,\dots,w_n$ may occur free, but $x,z,u$ do not occur free. Tarski arrives at this formula by considering the second-order continuity axiom:
\begin{equation}
\begin{split}
\forall XY(\exists z\forall xy(x\in X \land y\in Y \Rightarrow \be(zxy)) \\ 
\Rightarrow \exists u \forall xy (x\in X \land y \in Y \Rightarrow \be(xuy))).
\end{split}
\end{equation}

This axiom is second-order because it quantifies over sets. Tarski reformulates this as a schema for a (countably) infinite set of axioms, which as a whole mean the same thing as the second-order axiom, except \emph{only for sets definable in first-order logic}.

Given a real-closed ordered field $\mathbf{F}$, we can define the three-dimensional Euclidean geometry $\mathcal{E}_F=\mathbf{\langle F^3, \be_F, \ed_F\rangle}$, where $\mathbf{\be_F}$ is the 3-ary betweenness relation given by:
\begin{equation}
\mathbf{\be_F}=\left\{\langle \mathbf{a,b,c}\rangle: \begin{vmatrix} \mathbf{a_1} & \mathbf{a_2} & \mathbf{a_3} \\
\mathbf{b_1} & \mathbf{b_2} & \mathbf{b_3} \\
\mathbf{c_1} & \mathbf{c_2} & \mathbf{c_3} \end{vmatrix} = 0~\text{and}~( \mathbf{a_1\le b_1\le c_1}~\text{or}~\mathbf{c_1\le b_1\le a_1})\right\},
\end{equation}
and $\mathbf{\ed_F}$ is the 4-ary equidistance relation given by:
\begin{equation}
\begin{split}
\mathbf{\ed_F=\{\langle a,b,c,d\rangle:(a_1-b_1)^2+(a_2-b_2)^2+(a_3-b_3)^2} \\
\mathbf{=(c_1-d_1)^2+(c_2-d_2)^2+(c_3-d_3)^2\}.}
\end{split}
\end{equation}
In \cite{Tarski}, Tarski proves that for every real-closed ordered field $\mathbf{F}$, $\mathcal{E}_F$ is a model of $\mathfrak{E}$, and that every model $\M$ of $\mathfrak{E}$ is isomorphic to some $\mathcal{E}_F$. Since the ordered field $\mathcal{R}$ of real numbers is a real-closed ordered field, and every real-closed ordered field $F$ is elementary equivalent to $\mathcal{R}$, it follows that a sentence is true of any model of $\mathfrak{E}$ iff it is true of $\mathcal{E_R}$. Thus, $\mathfrak{E}$ is a complete theory and so by Proposition \ref{DecidabilityCriterion}, $\mathfrak{E}$ is decidable.



\section{Language and Definitions for \textsc{SimpleRel}}
\label{definitions}

We now define a first-order 2-sorted language $L_{SR}$ with variables of the sorts \emph{observer} and \emph{signal}, and containing two non-logical predicate symbols $T,R$ both of arity $2$. Observer variables will be written with lowercase Roman letters $a,b,c,d,\dots$, while signal variables will be written with lowercase greek letters $\alpha,\beta,\gamma,\dots$. We write $aT\alpha$ to express that an observer $a$ transmits a signal $\alpha$; and $bR\beta$ to express that an observer $b$ receives a signal $\beta$. In the canonical model, observers are interpreted as lines and signals are interpreted as line segments of slope $1$ (see section \ref{RCOFmodels}).

We define a number of relations in $L_{SR}$. Note that we are able to define the relations $\be$ and $\ed$ in terms of light signals. Also note that in many of these definitions, it is implicit that the speed of light is finite, and that it is in some sense invariant across reference frames. Pictures and short descriptions are given to help interpret these definitions. These are motivated by the interpretation of the definitions in the intended model $\M_F$, described later in section \ref{RCOFmodels}. Note also that some of these definitions, like the definition for $\sm$, are only meaningful for slower than light (STL) observers. These definitions and axioms are all that is needed to complete the proof that all observers are STL. These same comments also apply to the axioms given in the next section.

\begin{figure}
\begin{center}
\begin{tikzpicture}[scale=1]
\draw (0,0) node[right]{$a$} -- (.5,3);
\draw (-.5,0) node[left]{$d$} -- (2,3);
\draw (2,0) node[left]{$b$} -- (1.75,3);
\draw (2.5,0) node[right]{$c$} -- (0,2.75);
\end{tikzpicture}
\caption{$a$ and $b$ coplanar} \label{CopDef}
\end{center}
\end{figure}
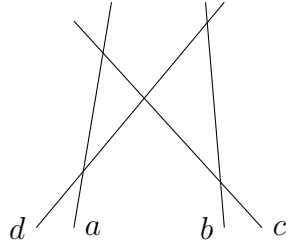

\medbreak
\noindent 
$\alpha$ is an \emph{event}. In the canonical model of \textsc{SimpleRel}, events are interpreted as points in spacetime:
\begin{equation}
\ev(\alpha)\equiv\forall a (aT\alpha \Leftrightarrow aR\alpha).
\end{equation}
$a$ \emph{meets} $b$:
\begin{equation}
{\rm{M}}(a,b)\equiv\exists\alpha(aT\alpha\land bT\alpha).
\end{equation}
$a$ and $b$ are \emph{coplanar} (Figure \ref{CopDef}):
\begin{equation}
\begin{split}
\cop(a,b)\equiv\exists c,d(c\ne d\land{\rm{M}}(a,c)\land{\rm{M}}(c,b)\land{\rm{M}}(d,b)\land{\rm{M}}(d,a) \\ 
\land\exists\gamma(cT\gamma\land dT\gamma\land\neg aT\gamma\land\neg bT\gamma))
\end{split}
\end{equation}

\noindent$a$ is \emph{parallel} to $b$. We later show in Lemma \ref{ParLemma} that parallelism is an equivalence relation:
\begin{equation}
a||b\equiv[\cop(a,b)\land\neg{\rm{M}}(a,b)]\lor[a=b].
\end{equation}
$\beta$ \emph{begins} at $\alpha$. In the next section we impose the axiom \textsc{AxEv}, from which we will have that every signal has a beginning and an end. Note that if $\alpha=\beg(\beta)$, then $\ev(\alpha)$:
\begin{equation}
\alpha=\beg(\beta)\equiv\ev(\alpha)\land \forall a(aT\beta\Leftrightarrow aT\alpha).
\end{equation}
$\beta$ \emph{ends} at $\alpha$. Note that if $\alpha=\nd(\beta)$, then $\ev(\alpha)$:
\begin{equation}
\alpha=\nd(\beta)\equiv\ev(\alpha)\land \forall a(aR\beta\Leftrightarrow aT\alpha).
\end{equation}
The segment $\alpha_1\alpha_2$ is \emph{lightlike}. Note that if $\rm{L}(\alpha_1,\alpha_2)$, $\ev(\alpha_1)$ and $\ev(\alpha_2)$:
\begin{equation}
\rm{L}(\alpha_1,\alpha_2)\equiv\exists\beta(\alpha_1=\beg(\beta)\land\alpha_2=\nd(\beta)).
\end{equation}
There is a light signal between $\alpha_1$ and $\alpha_2$:
\begin{equation}
l(\alpha_1,\alpha_2)\equiv \rm{L}(\alpha_1,\alpha_2)\lor\rm{L}(\alpha_2,\alpha_1)
\end{equation}
$b$ is \emph{between} $a$ and $c$ (Figure \ref{BeDef}):
\begin{equation}
\begin{split}
\be(abc)\equiv\exists\alpha,\beta,\gamma[a||b\land a||c
\land \rm{L}(\alpha,\beta)\land\rm{L}(\alpha,\gamma) \\ 
\land\rm{L}(\beta,\gamma)\land aT\alpha \land bT\beta\land cT\gamma].
\end{split}
\end{equation}

\begin{figure}
\begin{center}
\begin{tikzpicture}[scale=1]
\draw (0,0) node[left]{$a$} -- (0,3);
\draw (1,0) node[left]{$b$} -- (1,3);
\draw (2,0) node[left]{$c$} -- (2,3);
\fill (0,.5) circle (0.05) node[left]{$\alpha$};
\fill (1,1.5) circle (0.05) node[left]{$\beta$};
\fill (2,2.5) circle (0.05) node[left]{$\gamma$};
\draw[dashed] (0,.5) -- (2,2.5);
\end{tikzpicture}
\caption{$\be(abc)$} \label{BeDef}
\end{center}
\end{figure}

The distance between $a$ and $b$ is the same as the distance between $c$ and $d$ (Figure \ref{EqDef}):
\begin{equation}
\begin{split}
\ed(ab,cd)\equiv a||b\land a||c\land a||d \land \exists e\alpha_1\alpha_2\beta\gamma_1\gamma_2\delta\epsilon_1\epsilon_1'\epsilon_2\epsilon_2'[e||a \\ 
\land aT\alpha_1\land aT\alpha_2\land bT\beta\land cT\gamma_1\land cT\gamma_2\land dT\delta\land eT\epsilon_1 \land eT\epsilon_2 \\ 
\land eT\epsilon_1' \land eT\epsilon_2' \land L(\epsilon_1\alpha_1)\land L(\epsilon_1\gamma_1) \land L(\alpha_1\epsilon_1') \land L(\gamma_1\epsilon_1') \\
\land L(\epsilon_2\alpha_2)\land L(\epsilon_2\gamma_2) \land L(\alpha_2\epsilon_2') \land L(\gamma_2\epsilon_2')]
\end{split}
\end{equation}

\begin{figure}
\begin{center}
\begin{tikzpicture}[scale=1]
\draw (0,0) node[right]{$b$} -- (0,6);
\draw (1,0) node[right]{$a$} -- (1,6);
\draw (2.5,0) node[right]{$e$} -- (2.5,6);
\draw (4,0) node[right]{$c$} -- (4,6);
\draw (5,0) node[right]{$d$} -- (5,6);
\draw [dashed] (2.5,.5) node[right]{$\epsilon_1$} -- (1,2) node[left]{$\alpha_1$} -- (2.5,3.5) node[right]{$\epsilon_1'$} -- (4,2) node[right]{$\gamma_1$} -- (2.5,.5);
\draw [dashed] (1,2) -- (0,3) node[left]{$\beta$} -- (1,4) node[left]{$\alpha_2$} -- (2.5,2.5) node[right]{$\epsilon_2$} -- (4,4) node[right]{$\gamma_2$} -- (2.5,5.5) node[right]{$\epsilon_2'$} -- (1,4);
\draw [dashed] (4,4) -- (5,3) node[right]{$\delta$} -- (4,2);
\fill (2.5,.5) circle (0.05);
\fill (2.5,2.5) circle (0.05);
\fill (2.5,3.5) circle (0.05);
\fill (2.5,5.5) circle (0.05);
\fill (1,2) circle (0.05);
\fill (1,4) circle (0.05);
\fill (0,3) circle (0.05);
\fill (4,2) circle (0.05);
\fill (4,4) circle (0.05);
\fill (5,3) circle (0.05);
\end{tikzpicture}
\caption{$\ed(ab,cd)$} \label{EqDef}
\end{center}
\end{figure}

%

\noindent$\beta$ and $\beta'$ are \emph{simultaneous events} relative to $c$ (Figure \ref{SimDef}):
\begin{equation}
\begin{split}
\sm_c(\beta,\beta')\equiv\exists a \gamma_1\gamma_2\gamma_3\gamma_4[a||c\land aT\gamma_1\land\beg(\gamma_1)=\beg(\gamma_2) \\ 
\land \nd(\gamma_1)=\beg(\gamma_3)=\beta
\land\nd(\gamma_2)=\beg(\gamma_4)=\beta' 
\land aR\gamma_3 \\ \land \nd(\gamma_3)=\nd(\gamma_4)]
\end{split}
\end{equation}

\begin{figure}
\begin{center}
\begin{tikzpicture}[scale=1]
\draw (0,0) node[right]{$a$} -- (0,6);
\draw (5,0) node[right]{$c$} --  (5,6);
\draw[dashed] (0,1) -- node[right]{$\gamma_2$} (2,3) node[right] {$\beta'$} 
 -- node[right]{$\gamma_4$} (0,5) 
 -- node[left]{$\gamma_3$} (-2,3) node[left] {$\beta$}
 -- node[left]{$\gamma_1$} (0,1); 
\fill (2,3) circle (0.05);
\fill (-2,3) circle (0.05);
\end{tikzpicture}
\end{center}
\caption{$\sm_c(\beta,\beta')$} \label{SimDef}
\end{figure}

\begin{figure}
\begin{center}
\begin{tikzpicture}[scale=1]
\draw (0,0) node[right]{$b$} -- (0,5.5);
\draw (1,0) node[right]{$a$} -- (1,5.5);
\draw [dashed] (1,.5) node[right]{$\alpha_0'$} -- (0,1.5) node[left]{$\alpha_2$} -- (1,2.5) node[right]{$\alpha_1'$};
\draw [dashed] (1,3) node[right]{$\beta_0'$} -- (0,4) node[left]{$\beta_2$} -- (1,5) node[right]{$\beta_1'$};
\fill (3.5,.5) circle (.05) node[right]{$\alpha_0$};
\fill (4.5,2.5) circle (.05) node[right]{$\alpha_1$};
\fill (3,3) circle (0.05) node[right]{$\beta_0$};
\fill (5,5) circle (0.05) node[left]{$\beta_1$};
\fill (1,.5) circle (.05);
\fill (0,1.5) circle (.05);
\fill (1,2.5) circle (.05);
\fill (1,3) circle (.05);
\fill (0,4) circle (0.05);
\fill (1,5) circle (.05);
\end{tikzpicture}
\end{center}
\caption{$\Delta_a(\alpha_0\alpha_1,\beta_0\beta_1)$} \label{TimeDifDef}
\end{figure}
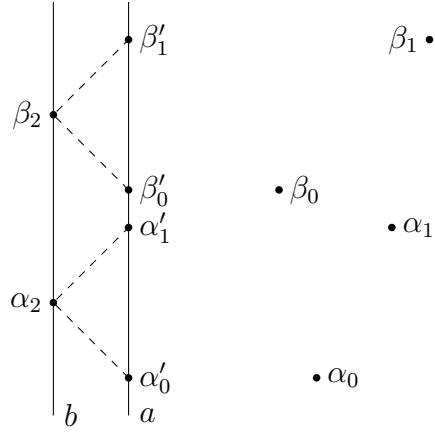

\noindent The difference in time between events $\alpha_0$ and $\alpha_1$ is the same as the difference in time between events $\beta_0$ and $\beta_1$, according to $a$ (Figure \ref{TimeDifDef}):
\begin{equation}
\begin{split}
\Delta_a(\alpha_0\alpha_1,\beta_0\beta_1)\equiv \exists b \alpha_0'\alpha_1'\beta_0'\beta_1'\{aT\alpha_0'\land aT\alpha_1'\land aT\beta_0'\land aT\beta_1' \\ 
\land \sm_a(\alpha_0,\alpha_0')\land\sm_a(\alpha_1,\alpha_1')\land\sm_a(\beta_0,\beta_0')\land\sm_a(\beta_1,\beta_1') \\
\land b||a \land [ (\alpha_0'=\alpha_1'\land\beta_0'=\beta_1')\lor(\alpha_0'\ne\alpha_1'\land\beta_0'\ne\beta_1'\land \\
\exists \alpha_2\beta_2(bT\alpha_2\land bT\beta_2\land l(\alpha_0',\alpha_2) \land l(\alpha_1',\alpha_2)\land l(\beta_0',\beta_2)\land l(\beta_1',\beta_2)))]\}
\end{split}
\end{equation}

$\alpha$ \emph{chronologically precedes} $\beta$, ie the interval from $\alpha$ to $\beta$ is timelike and $\alpha$ happens first. Note that if $\alpha\ll\beta$, then $\ev(\alpha)$ and $\ev(\beta)$:
\begin{equation}
\begin{split}
\alpha \ll\beta \equiv  \neg\rm{L}(\alpha,\beta)\land[\exists a,\gamma_1,\gamma_2(aT\alpha\land aT\beta \\
\land\alpha=\beg(\gamma_1)\land\nd(\gamma_1)=\beg(\gamma_2)\land\beta=\nd(\gamma_2))]
\end{split}
\end{equation}
$a$ is \emph{slower than light}. Note that relativity should have us say here: $a$ is slower than light according to some observer $b$. We can define the predicate $\rm{STL}$ in a global sense, though, by calling as STL those observers that receive a unique light signal from each event. In the canonical model, these are just those observers of slope greater than $1$ (Figure \ref{STLDef}):
\begin{equation}
\rm{STL}(a) \equiv \forall \gamma[\rm{Ev}(\gamma)\Rightarrow\exists!\beta(\beg(\beta)=\gamma\land aR\beta)]
\end{equation}

\begin{figure}
\begin{center}
\begin{tikzpicture}[scale=1]
\draw (0,0) -- (5,1) node[right]{$b$};
\draw [dashed](.833,.167) -- (2,-1) node[right]{$\gamma$} -- (3.75,.75);
\fill (2,-1) circle (0.05);
\end{tikzpicture}
\end{center}
\caption{$\neg\stl(b)$} \label{STLDef}
\end{figure}

\begin{figure}
\begin{center}
\begin{tikzpicture}[scale=1]
\draw (0,0) node[right]{$b$} -- (0,3);
\draw (1,0) node[right]{$a$} -- (1,3);
\draw (2,0) node[right]{$c$} -- (2,3);
\draw [dashed] (2,1) node[right]{$\gamma$} -- (1,2) node[left]{$\alpha_2$};
\fill (1,1) circle (.05) node[right]{$\alpha_1$};
\fill (2,1) circle (0.05);
\fill (1,2) circle (.05);
\end{tikzpicture}
\end{center}
\caption{$c=\tau_b(\alpha_1,\alpha_2)$} \label{TauDef}
\end{figure}
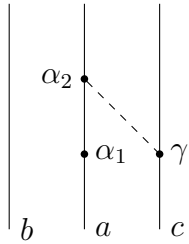

Finally, we give a technical definition of a symbol $\tau$ (Figure \ref{TauDef}). When $c=\tau_b(\alpha_1,\alpha_2)$, we consider the unique observer $a$ that transmits $\alpha_1$ and $\alpha_2$. Then $c$ is the unique observer that transmits some signal $\gamma$ such that $\nd(\gamma)=\alpha_2$ and $\sm_a(\beg(\gamma),\alpha_1)$, and such that $\be(bac)$. In the next section, \textsc{AxStIso} will give conditions under which $c$ exists. From Lemma \ref{AxLemma2}, we will see that $c$ is uniquely defined when it exists:
\begin{equation}
\begin{split}
c=\tau_b(\alpha_1,\alpha_2)\equiv\exists a \gamma[\stl(a)\land\be(bac)\land aT\alpha_1\land aT\alpha_2 \land \\
cT\gamma\land\sm_a(\beg(\gamma),\alpha_1)\land\nd(\gamma)=\alpha_2\land a\ne b\land\alpha_1\ll\alpha_2]
\end{split}
\end{equation}

\section{Axioms of \textsc{SimpleRel}}\label{STLaxioms}

In this section, we introduce \textsc{SimpleRel} as the following countably infinite system of axioms. See the previous section for some comments on interpretations of these axioms.

\textsc{Elementary~Geometry~Axioms (AxGeo)}: For each observer $a$, the observers moving parallel to $a$ satisfy the axioms $\mathfrak{E}$ of elementary geometry. For each of the axioms defined in section \ref{secEucGeo}, we add $\forall a$ to the beginning, then add, for each other observer mentioned in the axiom, that this observer is parallel to $a$. For instance, Axiom 1 becomes:
\begin{equation}
\forall a \forall xy[ x||a \land y||a \Rightarrow(\be(xyx)\Rightarrow(x=y))]
\end{equation}

\textsc{STL~observer~Axiom (AxSTL)}: There exists a slower than light observer. If $a$ is slower than light, and $b||a$, then $b$ is slower than light. 
\begin{equation} \exists a [\rm{STL}(a)]\land \forall a,b[\rm{STL}(a)\land b||a \Rightarrow \rm{STL}(b)]\end{equation}

\textsc{Events~Axiom (AxEv)}: Each signal is transmitted and received. The beginning and ends of signals are events.
\begin{equation}
\forall \alpha \exists a b\beta\gamma[aT\alpha \land bR\alpha\land\beg(\alpha)=\beta\land\nd(\alpha)=\gamma]
\end{equation}

\textsc{Time~Axiom (AxTime)}: Events transmitted by STL observers are ordered by $\ll$ and $=$. Events transmitted by non-STL observers are not relatable by $\ll$.
\begin{equation}
\begin{split}
\forall a \alpha\beta\{aT\alpha\land aT\beta\land \ev(\alpha)\land\ev(\beta)\Rightarrow[\rm{STL}(a) \\ 
\Rightarrow ((\alpha=\beta\land\neg\alpha\ll\beta \land\neg\beta\ll\alpha) \lor(\alpha\ll\beta\land\alpha\ne\beta\land\neg\beta\ll\alpha) \\ 
\lor(\beta\ll\alpha\land\alpha\ne\beta\land\neg\alpha\ll\beta)) \land \neg\rm{STL}(a)\Rightarrow(\neg\alpha\ll\beta\land\neg\beta\ll\alpha)]\}
\end{split}
\end{equation}

\textsc{Uniqueness~of~observers~Axiom (AxOb!)}: If $a$ and $b$ both transmit two signals $\alpha$ and $\beta$ such that $\beg(\alpha)\ne\beg(\beta)$, then $a=b$.
\begin{equation}
\forall a b \alpha \beta[\beg(\alpha)\ne\beg(\beta)\land aT\alpha\land aT\beta\land bT\alpha\land bT\beta\Rightarrow a=b]
\end{equation}

\textsc{Isotropy~Axiom (AxIso)}: Given an observer $a$ and an event $\beta$, there is a signal from $a$ to $\beta$ and a signal from $\beta$ to $a$.
\begin{equation}
\forall a\beta[\ev(\beta)\Rightarrow\exists\gamma_1\gamma_2(\beta=\nd(\gamma_1)\land aT\gamma_1\land\beta=\beg(\gamma_2)\land aR\gamma_2)]
\end{equation}

\textsc{Strong~Isotropy~Axiom (AxStIso)}: Given appropriate conditions, there exists $c$ such that $c=\tau_b(\alpha_1,\alpha_2)$.
\begin{equation}
\begin{split}
\forall a b\alpha_1\alpha_2[(aT\alpha_1\land aT\alpha_2\land\alpha_1\ll\alpha_2 \\
\land a\ne b \land a||b)\Rightarrow\exists c(c=\tau_b(\alpha_1,\alpha_2)]
\end{split}
\end{equation}

\begin{figure}
\begin{center}
\begin{tikzpicture}[scale=1]
\draw (0,0) node[right]{$b$} -- (0,3);
\draw (1,0) node[right]{$a$} -- (1,3);
\draw (2,0) node[right]{$c$} -- (2,3);
\draw (4,0) node[right]{$b'$} -- (4,3);
\draw (5,0) node[right]{$a'$} -- (5,3);
\draw (6,0) node[right]{$c'$} -- (6,3);
\draw [dashed] (2,1) -- (1,2) node[left]{$\alpha_2$};
\draw [dashed] (6,1) -- (5,2) node[left]{$\alpha_2'$};
\fill (2,1) circle (.05);
\fill (1,1) circle (0.05) node[left]{$\alpha_1$};
\fill (1,2) circle (0.05);
\fill (6,1) circle (.05);
\fill (5,1) circle (0.05) node[left]{$\alpha_1'$};
\fill (5,2) circle (0.05);
\end{tikzpicture}
\end{center}
\caption{\textsc{Position~Independence~of~Signal~Speed~Axiom (AxPoInd)}} \label{AxPoIndFig}
\end{figure}
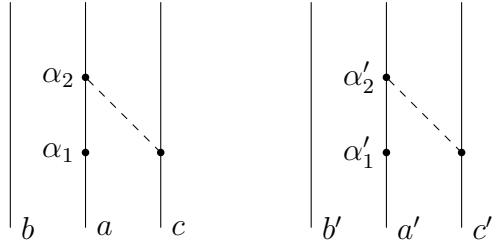

\textsc{Position~Independence~of~Signal~Speed~Axiom (AxPoInd)} (Figure \ref{AxPoIndFig}): The speed of light signals are independent of where and which direction they are sent in.
\begin{equation}
\begin{split}
\forall abca'b'c'\alpha_1\alpha_2\alpha_1'\alpha_2'\{[a||b\land a'||b' \land a||a' \land \sm_a(\alpha_1,\alpha_1') \\ \land aT\alpha_1 \land aT\alpha_2 \land a'T\alpha_1' \land a'T\alpha_2'\\
\land c=\tau_b(\alpha_1,\alpha_2)\land c'=\tau_{b'}(\alpha_1',\alpha_2')]\Rightarrow[\ed(ac,a'c')\Leftrightarrow\sm_a(\alpha_2,\alpha_2')]\} \\
\end{split}
\end{equation}

\textsc{Time~Independence~of~Signal~Speed~Axiom (AxTiInd)} (Figure \ref{AxTiIndFig}): The speed of light signals are independent of when they are sent.
\begin{equation}
\begin{split}
\forall acdd'\alpha_1\alpha_2\beta_1\beta_2\gamma_1\gamma_2\delta_1\delta_2\{[a||d\land a||d' \land d||c \land c\ne a \land \alpha_1=\beg(\gamma_1) \\
\land \nd(\gamma_1)=\beg(\gamma_2) \land cT\gamma_2 \land cT\delta_2 \land \beta_1=\beg(\delta_1) \land \nd(\delta_1)=\beg(\delta_2) \\ 
\land \alpha_2=\nd(\gamma_2) \land \beta_2=\nd(\delta_2)\land d=\tau_c(\alpha_1,\alpha_2) \land d'=\tau_{c}(\beta_1,\beta_2) \\ \land aT\alpha_1\land aT\alpha_2\land aT\beta_1\land aT\beta_2] \Rightarrow[d=d']\}
\end{split}
\end{equation}
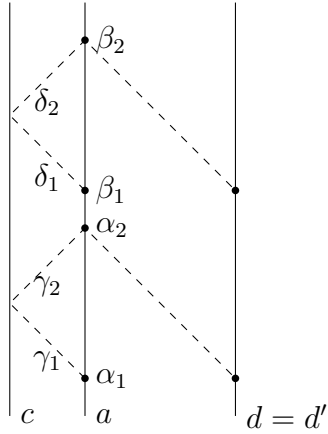
\begin{figure}
\begin{center}
\begin{tikzpicture}[scale=1]
\draw (0,0) node[right]{$c$} -- (0,5.5);
\draw (1,0) node[right]{$a$} -- (1,5.5);
\draw (3,0) node[right]{$d=d'$} -- (3,5.5);
\draw [dashed] (1,.5) node[right]{$\alpha_1$} -- node[below]{$\gamma_1$} (0,1.5) -- node[below]{$\gamma_2$} (1,2.5) node[right]{$\alpha_2$} -- (3,.5);
\draw [dashed] (1,3) node[right]{$\beta_1$} -- node [below]{$\delta_1$} (0,4) -- node[below]{$\delta_2$} (1,5) node[right]{$\beta_2$} -- (3,3);
\fill (1,.5) circle (0.05);
\fill (1,2.5) circle (0.05);
\fill (1,3) circle (0.05);
\fill (1,5) circle (0.05);
\fill (3,.5) circle (0.05);
\fill (3,3) circle (0.05);
\end{tikzpicture}
\end{center}
\caption{\textsc{Time~Independence~of~Signal~Speed~Axiom (AxTiInd)}} \label{AxTiIndFig}
\end{figure}

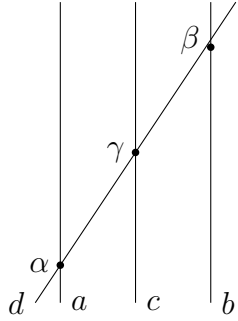
\begin{figure}
\begin{center}
\begin{tikzpicture}[scale=1]
\draw (0,0) node[right]{$a$} -- (0,4);
\draw (1,0) node[right]{$c$} -- (1,4);
\draw (2,0) node[right]{$b$} -- (2,4);
\draw (-.333,0) node[left]{$d$} -- (0,.5) node[left]{$\alpha$} -- (1,2) node[left]{$\gamma$} -- (2,3.5) node[left]{$\beta$} -- (2.333,4);
\fill (0,.5) circle (0.05);
\fill (1,2) circle (0.05);
\fill (2,3.4) circle (0.05);
\end{tikzpicture}
\end{center}
\caption{\textsc{Uniformity~of~observers~Axiom (AxUnOb)}} \label{AxUnObFig}
\end{figure}

\textsc{Uniformity~of~observers~Axiom (AxUnOb)} (Figure \ref{AxUnObFig}): For any two events $\alpha,\beta$ transmitted by an observer $d$, there is an event $\gamma$ transmitted by $d$ on the midpoint of the line between $\alpha$ and $\beta$. 
\begin{equation}
\begin{split}
\forall d \alpha \beta \{\ev(\alpha)\land\ev(\beta)\land dT\alpha\land dT\beta \Rightarrow \exists \gamma[\ev(\gamma) \land dT\gamma \land \forall abc (aT\alpha \\ 
\land bT\beta \land cT\gamma  \land a||b \land a||c \Rightarrow (\be(acb)\land \ed(ac,cb) 
 \land \Delta_a(\alpha\gamma,\gamma\beta)))]\}
\end{split}
\end{equation}

\textsc{Uniformity~of~Signal~Axiom (AxUnSi)}: Signals travel in straight lines.
\begin{equation}
\begin{split}
\forall abc\alpha\beta\gamma[\be(abc)\land aT\alpha\land bT\beta \\
\land cT\gamma\land L(\alpha,\gamma)\Rightarrow (L(\alpha,\beta)\Leftrightarrow L(\beta,\gamma))]
\end{split}
\end{equation}

\textsc{Rigid~Rods~Axiom (AxRR)}: Given a signal $\beta$ and an observer $a$, there exists an observer $b$ parallel to $a$ such that $bT\beta$. \begin{equation} \forall a \beta\exists b (bT\beta \land b||a)\end{equation}

\textsc{Simultaneity Axiom (AxSim)}: If $\sm_a(\alpha,\beta)$ and $\sm_a(\beta,\gamma)$, then $\sm_a(\alpha,\gamma)$.
\begin{equation}
\forall a \alpha\beta\gamma (\sm_a(\alpha,\beta)\land\sm_a(\beta,\gamma)\Rightarrow\sm_a(\alpha,\gamma))
\end{equation}

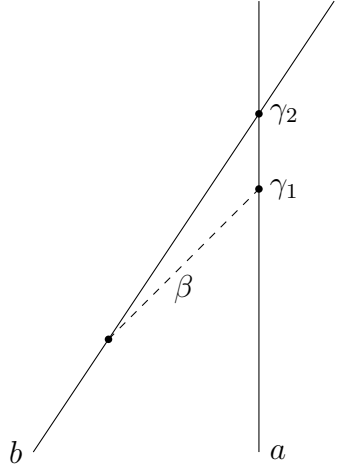
\begin{figure}
\begin{center}
\begin{tikzpicture}[scale=1]
\draw (-1,-1.5) node[left]{$b$} -- (2,3) node[right]{$\gamma_2$} -- (3,4.5);
\draw (2,-1.5) node[right]{$a$} -- (2,2) node[right]{$\gamma_1$} -- (2,4.5);
\draw [dashed] (0,0) -- node[below]{$\beta$} (2,2);
\fill (0,0) circle (0.05);
\fill (2,2) circle (0.05);
\fill (2,3) circle (0.05);
\end{tikzpicture}
\end{center}
\caption{\textsc{Limiting~Speed~Axiom (AxLim)}} \label{AxLimFig}
\end{figure}

\textsc{Limiting~Speed~Axiom (AxLim)}: Given an observer $a$ and signal $\beta$ such that $aR\beta$, for any event $\gamma_2$ transmitted by $a$ that occurs after $\nd(\beta)$, there is an observer $b$ that transmits $\beta$ and $\gamma_2$ (Figure \ref{AxLimFig}).
\begin{equation}
\begin{split}
\forall a \beta \gamma_1\gamma_2[(\gamma_1=\nd(\beta)\land\gamma_1\ll\gamma_2\land aT\gamma_1\land aT\gamma_2 ) \\ \Rightarrow\exists b (bT\beta\land bT\gamma_2)]
\end{split}
\end{equation}

\section{Useful Lemmas} \label{STLLemmas}

In this section, we prove a number of results that will be used to characterize all the models of \textsc{SimpleRel}. 

\begin{lemma} \label{ParLemma} Assume \textsc{AxGeo}. Then the relation $||$ is an equivalence relation between observers. \end{lemma}
\proof Reflexivity and symmetry of $||$ follow immediately from the definition. We show that if $a||b$ and $b||c$, then $a||c$. Since $b||a$ and $b||c$, it follows from \textsc{AxGeo} that $\ed(ac,ac)$ and hence by the definition of $\ed$, $a||c$.
\proofbox

\begin{coro} \label{AxLemma1} Assume \textsc{AxRR} and \textsc{AxGeo}. If $a$ is an observer and $\gamma$ is an event, there exists a \emph{unique} observer $b$ such that $bT\gamma$ and $b||a$.
\end{coro}
\proof $b$ exists by \textsc{AxRR}. Let $b'$ also be such that $b'||a$ and $b'T\gamma$. Then $b||b'$ by Lemma \ref{ParLemma}. Since $b$ meets $b'$ at $\gamma$, it follows from the definition of $||$ that $b=b'$. \proofbox

\begin{lemma} \label{UniquenessOfSignalsLemma} Assume \textsc{AxEv}, \textsc{AxSTL}, and \textsc{AxRR}. Then if two signals $\alpha$ and $\beta$ share endpoints, then $\alpha=\beta$. \end{lemma}
\proof Assume $\beg(\alpha)=\beg(\beta)$ and $\nd(\alpha)=\nd(\beta)$. By \textsc{AxEv}, let $\gamma=\nd(\alpha)=\nd(\beta)$. Consider the STL observer $b$ that transmits $\gamma$, which exists by \textsc{AxSTL} and \textsc{AxRR}. Since $b$ is STL, $\alpha=\beta$. \proofbox

\begin{lemma} \label{TwoEventsLemma} Assume \textsc{AxGeo}, \textsc{AxEv}, and \textsc{AxIso}. If $a$ is an observer, there are at least two distinct events $\alpha_1$ and $\alpha_2$ that are transmitted by $a$. \end{lemma}
\proof By Lemma \ref{ParLemma}, $a||a$ and hence, by a consequence of \textsc{AxGeo}, $\ed(aa,aa)$. By the definition of $\ed$, there is at least one event $\alpha_1$ transmitted by $a$. Assume that $\alpha_1$ is the only event transmitted by $a$. Then by the definition of $\ed$, \textsc{AxEv}, and \textsc{AxIso}, it follows that for every $b||a$, $\ed(ab,aa)$. But then $a=b$, by \textsc{AxGeo}, contradicting the dimension axioms of \textsc{AxGeo}. Thus, there is at least a second event $\alpha_2$ transmitted by $a$.
\proofbox

\begin{coro} \label{ObserversEventsLemma} Assume \textsc{AxGeo}, \textsc{AxEv}, \textsc{AxIso} and \textsc{AxOb!}. Then two observers $a$ and $b$ are identical iff for every signal $\alpha$, $aT\alpha\Leftrightarrow bT\alpha$. \end{coro}
\proof If $a=b$ then they transmit all the same signals. For observers $a$ and $b$, assume that given a signal $\alpha$, $aT\alpha\Leftrightarrow bT\alpha$. From Lemma \ref{TwoEventsLemma}, $a$ transmits at least two events. Thus, by \textsc{AxOb!}, $a=b$. \proofbox

\begin{lemma} \label{STLLemma1} Assume \textsc{AxIso}, \textsc{AxRR}, \textsc{AxUnSi}, \textsc{AxEv}, and \textsc{AxSTL}. Then if $a$ is STL, then for all events $\gamma$, there exists a unique signal $\beta$ such that $aT\beta$ and $\nd(\beta)=\gamma$ (Figure \ref{STLLemma1Fig}).
\end{lemma}

\proof By \textsc{AxIso}, such a signal exists. Assume that $aT\beta_1$ and $aT\beta_2$, and that $\nd(\beta_1)=\nd(\beta_2)=\gamma$. Let $c||a$ be such that $cT\gamma$, which exists by \textsc{AxRR}. Consider $b$ such that $\be(bac)$. Consider a signal $\alpha$ such that $bT\alpha$ and $\nd(\alpha)=\gamma$. Then, since $\rm{L}(\beg(\alpha),\gamma)$, $\rm{L}(\beg(\beta_1),\gamma)$, and $\rm{L}(\beg(\beta_2),\gamma)$, by \textsc{AxUnSi}, there exists signals $\alpha_1$ and $\alpha_2$ such that $\beg(\alpha_1)=\beg(\alpha_2)=\beg(\alpha)$, and $\nd(\alpha_1)=\beg(\beta_1)$ and $\nd(\alpha_2)=\beg(\beta_2)$. Since $a$ is STL, $\alpha_1=\alpha_2$ and hence by Lemma \ref{UniquenessOfSignalsLemma}, $\beta_1=\beta_2$. \proofbox

\begin{figure}
\begin{center}
\begin{tikzpicture}[scale=2]
\draw (0,0) node[right]{$b$} -- (0,3);
\draw (1,0) node[right]{$a$} -- (1,3);
\draw (2,0) node[right]{$c$} -- (2,3);
\draw [dashed] (0,.5) -- node[below]{$\alpha_1=\alpha_2$} (1,1.5) -- node[below]{$\beta_1=\beta_2$} (2,2.5) node[right]{$\gamma=\nd(\alpha)$};
\fill (0,.5) circle (0.025) node[left]{$\beg(\alpha)$};
\fill (1,1.5) circle (0.025);
\fill (2,2.5) circle (0.025);
\end{tikzpicture}
\end{center}
\caption{Proof of Lemma \ref{STLLemma1}} \label{STLLemma1Fig}
\end{figure}

\begin{lemma} \label{STLEventsLemma} Assume \textsc{AxIso}, \textsc{AxRR}, \textsc{AxUnSi}, \textsc{AxEv}, and \textsc{AxSTL}. If $a$ is STL, $aT\alpha_1$, $aT\alpha_2$, and $\sm_a(\alpha_1,\alpha_2)$, then $\alpha_1=\alpha_2$. \end{lemma}
\proof Since $\sm_a(\alpha_1,\alpha_2)$, there is an observer $b||a$ with signals $\beta_1,\beta_2$ such that $\beg(\beta_1)=\alpha_1$, $\beg(\beta_2)=\alpha_2$, and $\nd(\beta_1)=\nd(\beta_2)$. Let $\gamma=\nd(\beta_1)$. But since $a$ is STL, by Lemma \ref{STLLemma1}  we have that $\beg(\beta_1)=\beg(\beta_2)$ and hence $\alpha_1=\alpha_2$.
\proofbox

\begin{lemma} \label{AxLemma2} Assume \textsc{AxOb!}, \textsc{AxTime}, \textsc{AxPoInd}, and \textsc{AxGeo}. Then if $b=\tau_c(\alpha_1,\alpha_2)$ and $b'=\tau_c(\alpha_1,\alpha_2)$, then $b=b'$.
\end{lemma}
\proof Let $a$ be the observer that transmits $\alpha_1$ and $\alpha_2$. From \textsc{AxOb!}, $a$ is unique since $\alpha_1\ll\alpha_2$ and hence by \textsc{AxTime} $\alpha_1\ne\alpha_2$. Since $aT\alpha_1$ and $aT\alpha_2$, it follows from the definition of $\sm$ that $\sm_a(\alpha_1,\alpha_1)$ and $\sm_a(\alpha_2,\alpha_2)$. So, by \textsc{AxPoInd}, $\ed(ab,ab')$. Since $\be(cab)$ and $\be(cab')$, from \textsc{AxGeo} it follows that $b=b'$. \proofbox

\begin{lemma} \label{SimEventsLemma} Assume \textsc{AxRR}, \textsc{AxGeo}, \textsc{AxSTL}, \textsc{AxIso}, \textsc{AxUnSi}, \textsc{AxEv}, \textsc{AxOb!}, \textsc{AxTime}, \textsc{AxPoInd},  \textsc{AxTiInd}, and \textsc{AxSim}. For all events $\beta$ and all slower than light observers $a$ there exists a unique event $\alpha$ such that $aT\alpha$ and $\sm_a(\alpha,\beta)$ (Figure \ref{SimEventsLemmaFig}).
\end{lemma}
\proof Let $b$ be such that $b||a$ and $bT\beta$; $b$ exists by \textsc{AxRR}. Let $c$ be an observer such that $\ed(ca,cb)$, which exists by a consequence of \textsc{AxGeo}. Note that $b$ and $c$ are slower than light by \textsc{AxSTL}. Consider the unique (by Lemma \ref{STLLemma1}) signal $\gamma_1$ such that $cT\gamma_1$ and $\nd(\gamma_1)=\beta$; consider the unique signal $\gamma_2$ such that $\beg(\gamma_2)=\beg(\gamma_1)$ and $aR\gamma_2$; let $\alpha=\nd(\gamma_2)$. Consider $\gamma_3$ and $\gamma_4$ such that $\beg(\gamma_3)=\beta$, $\beg(\gamma_4)=\alpha$, $cR\gamma_3$, and $cR\gamma_4$. To show $\sm_a(\alpha,\beta)$, we show $\nd(\gamma_3)=\nd(\gamma_4)$.
 
Consider $d=\tau_a(\beg(\gamma_2),\nd(\gamma_4))$ and $e=\tau_b(\beg(\gamma_1),\nd(\gamma_3))$, which are uniquely defined by Lemma \ref{AxLemma2}. Considering the definition of $\ed$, since $\ed(ca,cb)$ there exists signals $\gamma_i'$, $i\in\{1,2,3,4\}$, such that $cT\gamma_1'$, $cT\gamma_2'$, $bR\gamma_1'$, $aR\gamma_2'$, $\nd(\gamma_1')=\beg(\gamma_3')$, $\nd(\gamma_2')=\beg(\gamma_4')$, $cR\gamma_3'$, and $cR\gamma_4'$. Furthermore, from the definition of $\ed$ there are signals $\gamma''_i$, $i\in\{1,2,3,4\}$ such that $\beg(\gamma'_i)=\beg(\gamma''_i)$ for $i=1,2$ and $\nd(\gamma'_i)=\beg(\gamma''_i)$ for $i=3,4$, and there is an observer $f||c$ such that $fR\gamma''_i$, $i\in\{1,2,3,4\}$. These $\gamma''_i$ are such that $\nd(\gamma''_1)=\nd(\gamma''_2)$ and $\nd(\gamma''_3)=\nd(\gamma''_4)$. But, by Lemma \ref{STLLemma1}, since $c$ is STL it follows that $\gamma''_1=\gamma''_2$ and $\gamma''_3=\gamma''_4$, and hence $\beg(\gamma_1')=\beg(\gamma_2')$, and $\nd(\gamma'_3)=\nd(\gamma'_4)$. Now, consider $d'=\tau_a(\beg(\gamma_2'),\nd(\gamma_4'))$ and $e'=\tau_b(\beg(\gamma_1'),\nd(\gamma_3'))$. By \textsc{AxPoInd}, $\ed(cd',ce')$. 

By \textsc{AxTiInd}, $d=d'$ and $e=e'$. By \textsc{AxPoInd}, $\sm_c(\nd(\gamma_3),\nd(\gamma_4))$ and hence $\nd(\gamma_3)=\nd(\gamma_4)$, because $c$ is STL. So $\sm_a(\alpha,\beta)$.

Now, assume that there exists an event $\alpha^*$ such that $aT\alpha^*$ and $\sm_a(\alpha^*,\beta)$. Then by \textsc{AxSim}, $\sm_a(\alpha,\alpha^*)$. By lemma \ref{STLEventsLemma}, $\alpha=\alpha^*$ and so $\alpha$ is unique.

\begin{figure}
\begin{center}
\begin{tikzpicture}[scale=1]
\draw (0,0) node[below]{$e=e'$} -- (0,5.5);
\draw (1,0) node[below]{$a$} -- (1,5.5);
\draw (2,0) node[below]{$c$} -- (2,5.5);
\draw (3,0) node[below]{$b$} -- (3,5.5);
\draw (4,0) node[below]{$d=d'$} -- (4,5.5);
\draw [dashed] (0,.5) -- (1,1.5) -- node[above]{$\gamma_4'$} (2,2.5) -- node[above]{$\gamma_3'$} (3,1.5) -- node[below]{$\gamma_1'$} (2,.5) -- node[below]{$\gamma_2'$} (1,1.5);
\draw [dashed] (4,.5) -- (3,1.5);
\draw [dashed] (0,3) -- (1,4) node[left]{$\alpha$} -- node[above]{$\gamma_4$} (2,5) -- node[above]{$\gamma_3$} (3,4) node[right]{$\beta$} -- node[below]{$\gamma_1$} (2,3) -- node[below]{$\gamma_2$} (1,4);
\draw [dashed] (4,3) -- (3,4);
\fill (0,.5) circle (0.05);
\fill (0,3) circle (0.05);
\fill (1,1.5) circle (0.05);
\fill (1,4) circle (0.05);
\fill (2,.5) circle (0.05);
\fill (2,2.5) circle (0.05);
\fill (2,3) circle (0.05);
\fill (2,5) circle (0.05);
\fill (3,1.5) circle (0.05);
\fill (3,4) circle (0.05);
\fill (4,.5) circle (0.05);
\fill (4,3) circle (0.05);
\end{tikzpicture}
\end{center}
\caption{Proof of Lemma \ref{SimEventsLemma}} \label{SimEventsLemmaFig}
\end{figure}

\section{Models over Real-Closed Ordered Fields}
\label{RCOFmodels}

In this section, we describe a class of structures which in the next section we show exhausts all the models of \textsc{SimpleRel}. Let $\mathbf{F}$ be a real-closed ordered field, and let $\lambda:\mathbf{F^4}\to\mathbf{F}$ be the Minkowski norm given by $\lambda(\mathbf{x_0,x_1,x_2,x_3})=\mathbf{-x_0^2+x_1^2+x_2^2+x_3^3}$. We say that a line $\mathbf{a}\subset\mathbf{F^4}$ is a \emph{timelike} line iff there exists $\mathbf{v,w\in F^4}$ such that $\lambda(\mathbf{w})\mathbf{<0}$ and $\mathbf{a=\{v+tw:t\in F\}}$. We say that $\mathbf{\alpha}$ is a \emph{lightlike interval} iff there exists $\mathbf{v,w\in F^4}$ and $\mathbf{t_1,t_2 \in F}$ with $\lambda(\mathbf{w})=\mathbf{0}$, $\mathbf{t_1\le t_2}$, and $\mathbf{\alpha=\{v+tw:t\in [t_1,t_2]\}}$. Let $\mathbf{\beg(\alpha)=v+t_1w}$ and $\mathbf{\nd(\alpha)=v+t_2w}$. Note that single points in $\mathbf{F^4}$ are lightlike intervals.

We now define a 2-sorted structure $\M_F$. Let the observers domain $\mathbf{Ob_F\subset \pow(F^4)}$ be the set of timelike lines of $\mathbf{F^4}$, and let the signals domain $\mathbf{Si_F\subset F^4\times F^4}$ be the set of lightlike intervals of $\mathbf{F^4}$. Let $\mathbf{T_F}$ be the 2-ary ``transmits'' relation consisting of the ordered pairs $\mathbf{\langle a, \alpha \rangle}$ where $\mathbf{a\in Ob_F}$, $\mathbf{\alpha\in Si_F}$, and $\mathbf{\beg(\alpha)\in a}$. Let $\mathbf{R_F}$ be the 2-ary ``receives'' relation consisting of the ordered pairs $\mathbf{\langle a, \alpha \rangle}$ $\mathbf{a\in Ob_F}$, $\mathbf{\alpha\in Si_F}$, and $\mathbf{\nd(\alpha)\in a}$. Let $\M_F=\mathbf{\langle Ob_F, Si_F,T_F,R_F\rangle}$.

In the Representation Theorem, we find that $\M_F\models \textsc{SimpleRel}$. Note that all observers in $\M_F$ are STL.

\section{Representation Theorem} \label{RepThmSTLSection}
\begin{theo} \label{RepThmSTL} $\M=\langle \mathbf{Ob, Si, T, R}\rangle$ is a model of \textsc{SimpleRel} iff there exists a real-closed ordered field $F$ such that $\M\cong\M_F$.
\end{theo}
\proof $\Leftarrow:$ One checks that $\M_\mathcal{R}$ satisfies each of the axioms in \textsc{SimpleRel}. Since $\mathbf{F}$ is elementary equivalent to $\mathcal{R}$, it follows that $\M_F$ is elementary equivalent to $\M_\mathcal{R}$, and hence $\M_F$ is a model of \textsc{SimpleRel}.

$\Rightarrow:$ By \textsc{AxSTL}, there exists a STL observer $\mathbf{a}$. Let $\mathbf{a}$ be fixed throughout this proof. Let $\mathbf{P_a}=\{\mathbf{b}\in\mathbf{P}:\mathbf{b||a}\}$.  We note that by \textsc{AxSTL}, for all $\mathbf{b\in P_a}$, $\mathbf{\rm{STL}(\mathbf{b})}$. Since \textsc{AxGeo} holds for $\mathbf{a}$, there exists a bijection $\phi:\mathbf{P_a}\to \mathbf{F^3}$, where $\mathbf{F^3}$ is a real-closed ordered field, inducing an isomorphism $\mathbf{\langle P_a; \be,\ed\rangle}\cong\mathbf{\langle F^3; \be_F, \ed_F\rangle}$. We construct an isomorphism between $\M$ and $\M_F$, which is a model of \textsc{SimpleRel}. 

We may assume without loss of generality that $\phi(\mathbf{a})=\mathbf{\langle 0,0,0\rangle}$, since for any $\mathbf{x,y\in F^3}$, there exists an automorphism of $\mathbf{\langle F^3; \be_F, \ed_F\rangle}$ taking $\mathbf{x}$ to $\mathbf{y}$. Let $\mathbf{\ev(\M)}=\{\mathbf{\alpha}:\M\mathcal{I}_{\alpha}^{\mathbf{\alpha}}~{\rm{sat}}~\ev(\alpha)\}.$ We define a bijection $\psi:\mathbf{\ev(\M)}\to \mathbf{F^4}$. Let $\mathbf{\theta}$ be an event such that $\mathbf{aT\theta}$ and let $\mathbf{\beta}$ be any event. $\theta$ exists by Lemma \ref{TwoEventsLemma}. We let $\mathbf{\theta}$ be fixed throughout this proof. Now, by \textsc{AxRR} there exists $\mathbf{b}\in\mathbf{P_a}$ such that $\mathbf{bT\beta}$. Corollary \ref{AxLemma1} gives us that $\mathbf{b}$ is unique. Since $\mathbf{a}$ is STL, by Lemma \ref{SimEventsLemma} there exists a unique event $\mathbf{\alpha}$ such that $\mathbf{aT\alpha}$ and $\mathbf{\sm_a(\alpha,\beta)}$. Let $\mathbf{\ev(a)}=\{\mathbf{\alpha}:\M\mathcal{I}_{a\alpha}^{\mathbf{a\alpha}}~\rm{sat}~\ev(\alpha)\land aT\alpha\}$. We define a function $\lambda:\mathbf{\ev(\M)}\to\mathbf{\ev(a)\times P_a}$ by $\lambda(\mathbf{\beta})=\mathbf{(\alpha,b)}$. $\lambda$ is injective by \textsc{TimeAx} and Lemma \ref{STLEventsLemma}. I claim that $\lambda$ is surjective. Indeed, given $\mathbf{(\alpha,b)}\in\mathbf{\ev(a)\times P_a}$, by Lemma \ref{SimEventsLemma} since $\mathbf{b}$ is STL there exists an event $\mathbf{\beta}$ such that $\mathbf{bT\beta}$ and $\mathbf{\sm_a(\alpha,\beta})$. And so $\lambda(\mathbf{\beta})=\mathbf{(\alpha,b)}$.

Let $\mathbf{a_{\pm1}}\in\mathbf{P_a}$ be the observer such that $\phi(\mathbf{a_{\pm1}})=\mathbf{\langle\pm1,0,0\rangle}\in\mathbf{F^3}$. Recalling that $\mathbf{aT\theta}$, we define $\mu:\mathbf{\ev(a)}\to\mathbf{P_a}$ by:
\begin{equation}
\mu(\mathbf{\alpha}) = \left\{
        \begin{array}{ll}
            \mathbf{\tau_{a_{-1}}(\theta,\alpha)} & \quad \rm{if}~\mathbf{\theta\ll\alpha} \\
            \mathbf{a} & \quad \rm{if}~\mathbf{\theta=\alpha} \\
            \mathbf{\tau_{a_{+1}}(\alpha,\theta)} & \quad \rm{if}~\mathbf{\alpha\ll\theta}
        \end{array}
    \right.
\end{equation}
Note that $\mu$ is well-defined for all $\mathbf{\alpha\in\ev(a)}$: by \textsc{AxTime}, since $\mathbf{a}$ is STL, $\mathbf{\ev(a)}$ is ordered by $\mathbf{\ll}$ and $=$, and by Lemma \ref{AxLemma2} $\mathbf{\tau}$ is a function. Let $\pi_1:\mathbf{F^3}\to\mathbf{F}$ be the projection of the first element, so $\pi_1(\mathbf{\langle e_1,e_2,e_3\rangle})=\mathbf{e_1}$. Define $\psi:\mathbf{\ev(\M)}\to\mathbf{F\times F^3}=\mathbf{F^4}$ so that the diagram in Figure \ref{ComDiagram} commutes.

\begin{figure}
\begin{center}
\begin{tikzpicture}
\node (D) at (0,2.5) {$\mathbf{Ev(\mathcal{M})}$};
\node (C) at (4.5,2.5) {$\mathbf{F\times F^3=F^4}$};
\node (B) at (4.5,0) {$\mathbf{F^3\times F^3}$};
\node (A) at (0,0) {$\mathbf{Ev(a)\times P_a}$};
\draw[->] (D) -- node[below]{$\psi$}  (C);
\draw[->] (D) -- node[right]{$\lambda$} (A);
\draw[->] (A) -- node[below]{$(\phi \circ \mu)\times \phi$} (B);
\draw[->] (B) -- node[right]{$\pi_1\times {\rm{Id}}$} (C);
\end{tikzpicture}
\end{center}
\caption{$\psi$ defined so the diagram commutes.} \label{ComDiagram}
\end{figure}

Note both functions $\mathbf{\tau_{a_{-1}}}$ and $\mathbf{\tau_{a_{+1}}}$ are injective because $\mathbf{a}$ is STL. We note that $\rm{RANGE}(\mathbf{\tau_{a_{-1}}})$ and $\rm{RANGE}(\mathbf{\tau_{a_{+1}}})$ are disjoint, because $\be(\mathbf{a_{-1},a,a_{+1}})$. Furthermore, $\mathbf{a}$ is not in either range, since otherwise $\theta=\alpha$. Thus, $\mu$ is injective. From the definition of $\mu$ we have that, for all $\mathbf{\alpha}\in\mathbf{\ev(a)}$, $\mathbf{\be(\mathbf{\mu(\alpha),a,a_{+1}})}$, $\mathbf{\be(\mathbf{\mu(\alpha),a_{+1},a})}$, or $\mathbf{\be(\mathbf{a,\mu(\alpha),a_{+1}})}$. Hence, $\pi_1\circ\phi\circ\mu$ is injective.

To show that $\psi$ is bijective, it remains to show that $\pi_1\circ\phi\circ\mu$ is surjective. This holds true if $\mathbf{F\times\{0\}\times\{0\}}\subset{\rm{RANGE}}(\phi\circ\mu)$. It then suffices to show that $\{\mathbf{b}: \M\mathcal{I}_{a,b,a_{+1}}^{\mathbf{a,b,a_{+1}}}~{\rm{sat}}~\be(b,a,a_{+1})\lor\be(a,b,a_{+1})\lor\be(a,a_{+1},b)\}\subset{\rm{RANGE}}(\mu)$. Indeed, consider $\mathbf{b}$ such that one of $\mathbf{b}$, $\mathbf{a}$, and $\mathbf{a_{+1}}$ is between the other two. Since $\mathbf{b}$ is STL, by Lemma \ref{SimEventsLemma}, there exists an event $\mathbf{\beta}$ transmitted by $\mathbf{b}$ such that $\mathbf{\sm_a(\theta,\beta)}$. Since $\mathbf{a}$ is STL, by \textsc{AxIso} there exists a unique signal $\mathbf{\gamma}$ from $\mathbf{\beta}$ to $\mathbf{a}$. Then $\mu(\mathbf{\nd(\gamma))}=\mathbf{b}$.

By Corollary \ref{ObserversEventsLemma}, \textsc{AxEv}, as well as Lemma \ref{UniquenessOfSignalsLemma}, we may identify observers in $\M$ and in $\M_F$ with the events they transmit, and we may identify signals with their endpoints. $\psi$ then induces functions $\psi_{Ob}:\mathbf{Ob \to \pow (F^4)}$ and $\psi_{Si}:\mathbf{Si \to F^4\times F^4}$. To show that $\M\cong\M_F$, we show that ${\rm{RANGE}}(\psi_{Ob})=\mathbf{Ob_F}$, ${\rm{RANGE}}(\psi_{Si})=\mathbf{Si_F}$, and that $\psi_{Ob}$ and $\psi_{Si}$ are bijections which preserve the relations $\mathbf{T,R}$.

Now, let $\mathbf{b}\in\mathbf{Ob}$ such that $\psi_{Ob}(\mathbf{b})\in \mathbf{Ob_F}$, and let $\mathbf{\beta}\in\mathbf{Si}$ such that $\psi_{Si}(\mathbf{\beta})\in \mathbf{Si_F}$. Let $\mathbf{\gamma=\beg(\beta)}$. From this identification, $\mathbf{bT\beta}\Leftrightarrow\mathbf{bT\gamma}$. Since $\gamma$ is an event, this is true iff $\psi_{Ob}(\mathbf{b})\mathbf{T_F}\psi_{Si}(\mathbf{\gamma})$ or equivalently $\psi_{Ob}(\mathbf{b})\mathbf{T_F}\psi_{Si}(\mathbf{\beta})$. Thus, $\mathbf{bT\beta}\Leftrightarrow\mathbf{\psi_{Ob}(b)T_F\psi_{Si}(\beta)}$. A similar argument shows that $\mathbf{bR\beta}\Leftrightarrow\mathbf{\psi_{Ob}(b)R_F\psi_{Si}(\beta)}$. This fact will be used several times in special cases in the remainder of this proof. Additionally, it gives us that to show that $\psi_{Ob}$ and $\psi_{Si}$ preserve the relations $\mathbf{T,R}$, it suffices to show that ${\rm{RANGE}}(\psi_{Ob})=\mathbf{Ob_F}$, ${\rm{RANGE}}(\psi_{Si})=\mathbf{Si_F}$, and that $\psi_{Ob}$ and $\psi_{Si}$ are bijections.

Let $\mathbf{\gamma}\in \mathbf{Si}$. Then there exists events $\mathbf{\alpha=\beg(\gamma)}$ and $\mathbf{\beta=\nd(\gamma)}$ such that $\mathbf{\alpha\beta}$ is lightlike. Let $\psi(\mathbf{\alpha})=\mathbf{\langle a_0,a_1,a_2,a_3\rangle}$ and $\psi(\mathbf{\beta})=\mathbf{\langle b_0,b_1,b_2,b_3\rangle}$ (Figure \ref{SigRepThmFig}). We give the proof in the case that $a_0,b_0>0$; the other cases are very similar.
By the definition of $\psi$, as well as \textsc{AxSim}, any two events $\mathbf{\epsilon,\epsilon'}$ with $\psi(\mathbf{\epsilon})=\mathbf{\langle e_0,e_1,e_2,e_3\rangle}$ and $\psi(\mathbf{\epsilon'})=\mathbf{\langle e'_0,e'_1,e'_2,e'_3\rangle}$ are such that $\sm_a(\epsilon,\epsilon')$ iff $\mathbf{e_0=e'_0}$. Let $\mathbf{\beta',\gamma,\gamma'}$ be such that  $\psi(\mathbf{\beta'})=\mathbf{\langle b_0,0,0,0\rangle}$, $\psi(\mathbf{\gamma})=\mathbf{\langle a_0,b_1,b_2,b_3\rangle}$ and $\psi(\mathbf{\gamma'})=\mathbf{\langle a_0,0,0,0\rangle}$. Thus, by \textsc{AxPoInd} as well as how $\psi$ is defined, since $\mathbf{\sm_a(\beta,\beta')}$ and $\mathbf{\sm_a(\gamma,\gamma')}$, we have that $\mathbf{\alpha'\beta'}$ is lightlike, where $\psi(\alpha')=\mathbf{\langle a_0,x,0,0\rangle}$ and $\mathbf{x=\sqrt{(a_1-b_1)^2+(a_2-b_2)^2+(a_3-b_3)^2}}$. Then, considering the definition of $\psi$ gives us that $\mathbf{L(\psi^{-1}(\langle 0,b_0,0,0\rangle),\beta')}$ and $\mathbf{L(\psi^{-1}(\langle 0,a_0,0,0\rangle),\gamma')}$. And, by \textsc{AxUnSi} $\mathbf{L(\psi^{-1}(\langle 0,b_0,0,0\rangle),\alpha')}$. But then since $\mathbf{\sm_a(\psi^{-1}(\langle 0,0,0,0\rangle),\psi^{-1}(\langle 0,x,0,0\rangle)}$ and $\mathbf{\sm_a(\gamma',\alpha')}$, we have by \textsc{AxPoInd} that $\mathbf{a_0=b_0-x}$. It follows that $\mathbf{a_0\le b_0}$ and $\mathbf{a_0=b_0}$ iff $\mathbf{x=0}$ and $\mathbf{\alpha=\beta}$. It also follows that, 
\begin{equation}
\mathbf{(b_0-a_0)^2=(b_1-a_1)^2+(b_2-a_2)^2+(b_3-a_3)^2},
\end{equation}
and hence the line segment defined by $\psi(\mathbf{\alpha})\psi(\mathbf{\beta})$ is lightlike in $\mathbf{F^4}$. Thus, $\psi_{Si}(\gamma)\in\mathbf{Si_F}$, and so ${\rm{RANGE}}(\psi_{Si})=\mathbf{Si_F}$. And $\psi_{Si}$ is one-to-one because $\psi$ is one-to-one.

\begin{figure}
\begin{center}
\begin{tikzpicture}[scale=1.5]
\draw [dashed] (8,5) node[right]{$\psi(\alpha)$} -- (6,7) node[right]{$\psi(\beta)$};
\fill (8,5) circle (0.05);
\fill (6,7) circle (0.05);
\fill (6,5) circle (0.05) node[right]{$\psi(\gamma)$};
\draw [dashed] (7,0) node[below right]{$\mathbf{\langle 0,b_0,0,0\rangle}$} -- (2,5) node[right]{$\psi(\alpha')$} -- (0,7) node[right]{$\psi(\beta')$};
\fill (7,0) circle (0.05);
\fill (0,7) circle (0.05);
\fill (2,5) circle (0.05);
\draw [dashed] (5,0) node[below]{$\mathbf{\langle 0,a_0,0,0\rangle}$} -- (0,5) node[right]{$\psi(\gamma')$};
\fill (0,5) circle (0.05);
\fill (5,0) circle (0.05);
\fill (2,0) circle (0.05) node[below]{$\mathbf{\langle 0,x,0,0 \rangle}$};
\fill (0,0) circle (0.05) node[below]{$\mathbf{\langle 0,0,0,0\rangle}$};
\draw [dotted] (7,0) -- (0,0) -- (0,7) -- (6,7);
\draw [dotted] (0,5) -- (8,5);
\draw [dotted] (2,0) -- (2,5);
\end{tikzpicture}
\end{center}
\caption{Argument that, if $\mathbf{\alpha\beta}$ is lightlike in $\M$, then $\psi(\mathbf{\alpha})\psi(\mathbf{\beta})$ is lightlike in $\M_F$} \label{SigRepThmFig}
\end{figure}

Now, let $\mathbf{\gamma'\in Si_F}$, and let its endpoints be the events $\psi(\mathbf{\alpha})$ and $\psi(\mathbf{\beta})$. We show that in $\M$, $\mathbf{\rm{L}(\alpha,\beta)}$. Let $\mathbf{b}\in\mathbf{P_a}$ be such that $\mathbf{bT\beta}$. Note that by construction of $\psi$, $\psi(\ev(\mathbf{b}))$ is a line in $\mathbf{F^4}$ given by $\mathbf{\langle t,x_b,y_b,z_b\rangle}$ where $\mathbf{t}$ ranges over all values of $\mathbf{F}$ and $\mathbf{x_b,y_b,z_b}$ are fixed. Thus, $\psi_{Ob}(\mathbf{b})\in\mathbf{Ob_F}$ is STL. Now, by \textsc{AxIso}, there exists a signal $\mathbf{\gamma}$ from $\mathbf{\alpha}$ to $\mathbf{b}$. We have that $\psi(\mathbf{\alpha})\psi(\mathbf{\nd(\gamma)})$ is lightlike, so since $\M_F\models \textsc{SimpleRel}$, and since $\psi_{Ob}(\mathbf{b})$ is STL, we must have $\psi_{Si}(\mathbf{\gamma})=\mathbf{\gamma'}$. Thus, $\psi_{Si}$ is onto. 

We now turn to observers. Let $\mathbf{b\in Ob}$ that transmits events $\mathbf{\alpha_0,\alpha_1,\alpha_2}$. I claim that $\psi(\mathbf{\alpha_0}),\psi(\mathbf{\alpha_1}),\psi(\mathbf{\alpha_2})$ are collinear. Assume they are not. Using the fact that for all $\mathbf{c\in P_a}$, $\psi_{Ob}(\mathbf{c})\in\mathbf{Ob_F}$ and $\psi(\ev(\mathbf{c}))$ is a vertical line in $\mathbf{F^4}$, by \textsc{AxUnOb}, we find that for any two events $\mathbf{\alpha,\alpha'}$ transmitted by $\mathbf{b}$, there is an event $\mathbf{\alpha''}$ with $\mathbf{bT\alpha''}$ and such that $\psi(\mathbf{\alpha''})$ lies midway on the line segment formed by $\psi(\mathbf{\alpha})$ and $\psi(\mathbf{\alpha'})$. By applying this fact, we may find events transmitted by $\psi(\mathbf{b})$ that evenly subdivide each of the sides of the triangle formed by $\psi(\mathbf{\alpha_0}),\psi(\mathbf{\alpha_1}),\psi(\mathbf{\alpha_2})$ into $2^n$ segments, for any $n\in\N$. 

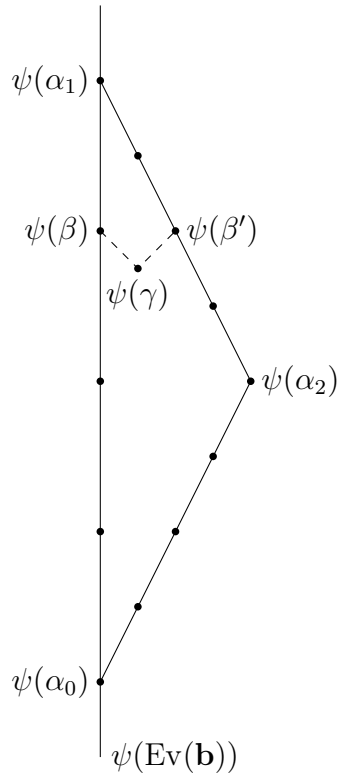
\begin{figure}
\begin{center}
\begin{tikzpicture}[scale=2]
\draw (0,-.5) node[right]{$\psi(\ev(\mathbf{b}))$} -- (0,0) node[left]{$\psi(\alpha_0)$}  -- (0,3) node[left]{$\psi(\beta)$} -- (0,4) node[left]{$\psi(\alpha_1)$} -- (.5,3) node[right]{$\psi(\beta')$} -- (1,2) node[right]{$\psi(\alpha_2)$} -- (0,0);
\draw (0,4) -- (0,4.5);
\draw [dashed] (0,3) -- (.25,2.75) node[below]{$\psi(\gamma)$} -- (.5,3);
\fill (0,0) circle (.025);
\fill (0,1) circle (.025);
\fill (0,2) circle (.025);
\fill (0,3) circle (.025);
\fill (0,4) circle (.025);
\fill (.5,3) circle (.025);
\fill (.25,.5) circle (.025);
\fill (.5,1) circle (.025);
\fill (.75,1.5) circle (.025);
\fill (.75,2.5) circle (.025);
\fill (.25,3.5) circle (.025);
\fill (.25,2.75) circle (.025);
\fill (1,2) circle (.025);
\end{tikzpicture}
\end{center}
\caption{Construction of $\psi(\beta)$ and $\psi(\beta')$ in the STL case, by subdividing the triangle formed by $\psi(\mathbf{\alpha_0}),\psi(\mathbf{\alpha_1}),\psi(\mathbf{\alpha_2})$ with the images of events transmitted by $\mathbf{b}$.} \label{ObConstructionFig1}
\end{figure}

If $\mathbf{b}$ is STL, we can subdivide the sides of this triangle until we find two events $\mathbf{\beta\ne\beta'}$ transmitted by $\mathbf{b}$ such that there is an event $\psi(\mathbf{\gamma})\in\mathbf{Si_F}$ with $\mathbf{\rm{L}(\psi(\mathbf{\gamma}),\psi(\mathbf{\beta}))}$ and $\mathbf{\rm{L}(\psi(\mathbf{\gamma}),\psi(\mathbf{\beta'}))}$ (Figure \ref{ObConstructionFig1}). Since $\psi_{Si}$ is a bijection, $\mathbf{\rm{L}(\mathbf{\gamma},\mathbf{\beta})}$ and $\mathbf{\rm{L}(\mathbf{\gamma},\mathbf{\beta'})}$ in $\M$, contradicting the fact that $\mathbf{b}$ is STL.

If $\mathbf{b}$ is not STL, the argument goes the same way, except we find $\beta$ and $\beta'$ such that $\mathbf{\beta\ll\beta'}$, contradicting \textsc{AxTime}. Thus, there exists a line $l(\mathbf{b})$ in $\mathbf{F^4}$ such that $\psi(\ev(\mathbf{b}))\subset l(\mathbf{b})$. 


We now see that, if for an observer $\mathbf{b}$, $\psi_{Ob}(\mathbf{b})$ is contained within a non-timelike line in $\mathbf{F^4}$, then \textsc{AxIso} fails to hold for $\mathbf{b}$, since $\psi_{Si}$ is a bijection. Thus, all observers in $\M$ are STL. Note that this result obtains because of assumptions about the isotropy of spacetime. Also, it follows from \textsc{AxIso}, as well as the fact that the image of any observer lies on a timelike line, for an observer $\mathbf{b}$, if $\psi(\ev(\mathbf{b}))\subset l(\mathbf{b})$ in $\mathbf{F^4}$, then $\psi_{Ob}(\mathbf{b})=l(\mathbf{b})$. So for all $\mathbf{b}\in\M$, $\psi_{Ob}(\mathbf{b})$ is an observer in $\M_F$. Hence, ${\rm{RANGE}}(\psi_{Ob})=\mathbf{Ob_F}$. And $\psi_{Ob}$ is one-to-one because $\psi$ is one-to-one.

Let $\mathbf{b}$ be a (STL) observer in $\M_F$, and let $\psi(\mathbf{\mu_1})\ll\psi(\mathbf{\mu_2})$ be two events that $\mathbf{b}$ transmits. We will show there exists an observer $\mathbf{b'\in Ob}$ such that $\psi_{Ob}(\mathbf{b'})=\mathbf{b}$. Consider the line $\psi_{Ob}(\mathbf{c})$ where $\mathbf{c\in P_a}$ and $\mathbf{cT\mu_2}$. By the structure of $\M_F$, there exists an event $\psi(\mathbf{\gamma})$ such that $\psi_{Ob}\mathbf{(\mathbf{c})T\psi(\gamma)}$, $\psi(\mathbf{\gamma})\mathbf{\ll}\psi(\mathbf{\mu_2})$, and $\mathbf{\rm{L}}(\psi(\mathbf{\mu_1}),\psi(\mathbf{\gamma}))$. Since $\psi_{Si}$ is a bijection, by \textsc{AxLim}, there is an observer $\mathbf{b'}$ in $\M$ such that $\mathbf{b'T\mu_1}$ and $\mathbf{b'T\mu_2}$. By the previous paragraph and \textsc{AxOb!}, we have $\psi_{Ob}(\mathbf{b'})=\mathbf{b}$ and hence $\psi_{Ob}$ is onto.  Thus, $\psi_{Ob}$ is a bijection.

As was noted before, showing $\psi_{Ob}:\mathbf{Ob\to Ob_F}$ and $\psi_{Si}:\mathbf{Si\to Si_F}$ are bijections suffices to show that $\psi_{Ob}$ and $\psi_{Si}$ preserve the relations $\mathbf{T,R}$. Hence, $\M\cong\M_F$. \proofbox

\begin{coro} It is provable from \textsc{SimpleRel} that all observers are slower than light. \end{coro}
\proof Follows from the Representation Theorem  and the Completeness Theorem (Theorem \ref{ComTheo}). \proofbox

\begin{coro} \label{Coro2STL} The set of consequences of \textsc{SimpleRel} is complete and decidable. \end{coro}
\proof From the Representation Theorem, the set of consequences of \textsc{SimpleRel} is the theory of the structure $\M_{\mathcal{R}}$. Since, for any sentence $\Phi$ in the language of \textsc{SimpleRel}, $\M_{\mathcal{R}}\models\Phi$ or $\M_{\mathcal{R}}\models\neg\Phi$, the consequences of \textsc{SimpleRel} are complete. By Proposition \ref{DecidabilityCriterion}, since the consequences of \textsc{SimpleRel} are axiomatized by \textsc{SimpleRel}, they are decidable. \proofbox

\chapter{A simple axiomatization of special relativity with faster-than-light observers} \label{FTLchapter}

\section{Introduction}
In the last section, we introduced \textsc{SimpleRel}, an axiomatization for special relativity without faster than light observers. In this section, we develop \textsc{SimpleRelFTL}, a modification of \textsc{SimpleRel} that allows for the existence of faster than light observers. We prove a representation theorem analogous to the one for \textsc{SimpleRel}. We make only a few, relatively small, modifications to \textsc{SimpleRel} to arrive at \textsc{SimpleRelFTL}. The proof of the representation theorem for \textsc{SimpleRelFTL} is also fairly similar to the STL case, as most of it relies on the existence of just STL observers.

\section{Definitions for \textsc{SimpleRelFTL}} \label{definitionsFTL}
Like we did in section \ref{definitions}, we now give a language and set of definitions that will be used in the axiom system \textsc{SimpleRelFTL}. We use the same first-order $2$-sorted language as before. Some definitions will be changed here in order to make them work for FTL observers, and some new ones will be added. The defined terms \emph{event, meets, coplanar, parallel, begins, ends, lightlike, existence of a light signal, chonologically precedes}, and \emph{STL} remain the same as before. The predicates \emph{between, equidistance, simultaneous, equidistance in time,} and $\tau$ are given different definitions. We denote the newly defined predicates $\be^{FTL},\ed^{FTL},\sm^{FTL} \Delta^{FTL}$, and $\tau^{FTL}$. The chief difficulty here is that not all observers can send light signals to one another. In section \ref{FTLLemmas}, we prove that the new definitions are equivalent to the original definitions for STL observers in \textsc{SimpleRelFTL}. We also note in section \ref{RepThmFTL} that we could have used the changed definitions as our definitions for terms in \textsc{SimpleRel}. 

It is worth noting the motivation behind these changed definitions. Inspection of the proof of theorem \ref{RepThmFTL} reveals that the only times in which we use the relations $\be^{FTL},\ed^{FTL},\sm^{FTL} \Delta^{FTL}$, and $\tau^{FTL}$ is with respect to STL observers. Having a meaningful definition of these terms for FTL observers is not, then, necessary to establish a representation theorem for \textsc{SimpleRelFTL}. This is reflected by the fact that we find that the relatable dual $(a)_b'$ of $a$ with respect to $b$ is well-defined only as a corollary of theorem \ref{RepThmFTL}; this fact was not needed at any point in the proof of the representation theorem. The reason we revise definitions here is, then, not to ensure that the theory works as it should. Rather, we revise definitions here so that the definitions reflect what we actually mean by betweenness, Minkowskian equidistance, etc., and do so for \emph{all} observers, not just STL ones.

\medbreak
\noindent $a$ is \emph{lightspeed}:
\begin{equation}
\exists \alpha (\neg\ev(\alpha)\land aT\alpha\land aR\alpha)
\end{equation}

\noindent $a$ is \emph{faster than light}:
\begin{equation}
\ftl(a)\equiv \neg\lightspeed(a)\land\neg\stl(a)
\end{equation}

\noindent $a$ \emph{transmits or receives} $\alpha$
\begin{equation}
aTR\alpha\equiv aT\alpha\lor aR\alpha
\end{equation}

\noindent $a$ and $b$ are \emph{relatable} by light signals:
\begin{equation}
\rho(a,b)\equiv\exists \gamma(aTR\gamma\land bTR\gamma)
\end{equation}

\begin{figure}
\begin{center}
\begin{tikzpicture}[scale=2]
\draw (1,0) node[left]{$b$} -- (3,1);
\draw (0,1) node[left]{$a$} -- (2,2);
\draw [dashed] (2,.5) -- (1,1.5);
\draw (2,1.5) ellipse (1 and .1);
\draw[dashed](2,.5) -- (3,1.5);
\draw [dashed](2,.5) -- (1,-.5);
\draw [dashed] (2,.5) -- (3,-.5);
\draw (2,-.5) ellipse (1 and .1);
\end{tikzpicture}
\end{center}
\caption{$OP(ab)$. The diagram is drawn in three dimensions, and the cone represents a light cone emitted by $b$. The plane formed by $a$ and $b$ is tangent to all the light cones emitted by $a$ and $b$.} \label{OPFig}
\end{figure}
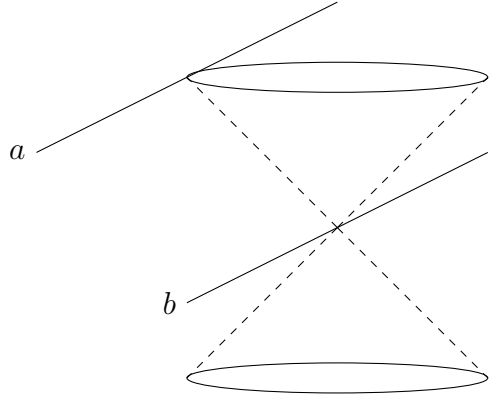

\noindent The plane defined by $a$ and $b$ is tangent to their light cones. It is called an \emph{optical plane} (Figure \ref{OPFig}):
\begin{equation}
OP(a,b)\equiv a||b\land \exists \gamma(aTR\gamma\land bTR\gamma)\land \forall c(M(a,c)\land M(b,c)\Rightarrow \ftl(c))
\end{equation}

\noindent $a,b$ and $c$ are relatable, and $b$ is between $a$ and $c$ (Figure \ref{BeDef})
\begin{equation}
\begin{split}
\be_\rho(abc)\equiv a||b\land a||c \land \exists \alpha\beta\gamma[aT\alpha\land bT\beta\land cT\gamma  \\ 
\land((L(\alpha\beta)\land L(\alpha\gamma)\land L(\beta\gamma))\lor(L(\gamma\beta)\land L(\gamma\alpha)\land L(\beta\alpha)))]
\end{split}
\end{equation}

\noindent $a$ and $b$ are relatable, and $c$ and $d$ are relatable, and the distance between $a$ and $b$ is the same as the distance between $c$ and $d$. There are three cases here. First, if $a=b$ and $c=d$, then $\ed_\rho(ab,cd)$ holds. Second, if the pairs $a,b$ and $c,d$ each define an optical plane, then $\ed_\rho(ab,cd)$ holds. This is motivated by the fact that, in Minkowskian spacetime, lightlike intervals all have length $0$. Third, $\ed_\rho(ab,cd)$ holds if we can find some observer $e||a$ such that there are light signals between $e$ and both $a$ and $c$, and such that we can send light signals between $a,b,c,d,e$ so that they always meet in a certain way. (Figure \ref{EqDefFTL}):
\begin{equation}
\begin{split}
\ed_\rho(ab,cd)\equiv a||b\land a||c \land a||d \\ 
\land \{[a=b\land c=d]\lor[OP(ab)\land OP(cd)] \\
\lor \exists \alpha_1\alpha_2\beta\gamma_1\gamma_2\delta\{\alpha_1\ne\alpha_2\land\gamma_1\ne\gamma_2\land l(\alpha_1\beta)\land l(\alpha_2\beta)\land l(\gamma_1\delta)\land l(\gamma_2\delta) \\ \land aT\alpha_1\land aT\alpha_2
\land bT\beta \land cT\gamma_1\land cT\gamma_2\land dT\delta \\ \land \exists e[e||a\land\exists \epsilon(eT\epsilon\land l(\epsilon\alpha_1)\land l(\epsilon\gamma_1)) 
\\ \land \forall \epsilon[eT\epsilon \Rightarrow ((l(\epsilon\alpha_1)\Leftrightarrow l(\epsilon\gamma_1))\land(l(\epsilon\alpha_2)\Leftrightarrow l(\epsilon\gamma_2)))]]\}\}
\end{split}
\end{equation}

\begin{figure}
\begin{center}
\begin{tikzpicture}[scale=1]
\draw (0,0) node[right]{$b$} -- (0,6);
\draw (1,0) node[right]{$a$} -- (1,6);
\draw (2.5,0) node[right]{$e$} -- (2.5,6);
\draw (4,0) node[right]{$c$} -- (4,6);
\draw (5,0) node[right]{$d$} -- (5,6);
\draw [dashed] (2.5,.5) -- (1,2) node[left]{$\alpha_1$} -- (2.5,3.5) -- (4,2) node[right]{$\gamma_1$} -- (2.5,.5);
\draw [dashed] (1,2) -- (0,3) node[left]{$\beta$} -- (1,4) node[left]{$\alpha_2$} -- (2.5,2.5) -- (4,4) node[right]{$\gamma_2$} -- (2.5,5.5) -- (1,4);
\draw [dashed] (4,4) -- (5,3) node[right]{$\delta$} -- (4,2);
\fill (2.5,.5) circle (0.05);
\fill (2.5,2.5) circle (0.05);
\fill (2.5,3.5) circle (0.05);
\fill (2.5,5.5) circle (0.05);
\fill (1,2) circle (0.05);
\fill (1,4) circle (0.05);
\fill (0,3) circle (0.05);
\fill (4,2) circle (0.05);
\fill (4,4) circle (0.05);
\fill (5,3) circle (0.05);
\end{tikzpicture}
\caption{The third case of $\ed(ab,cd)$.} \label{EqDefFTL}
\end{center}
\end{figure}

\begin{figure}
\begin{center}
\begin{tikzpicture}[scale=.8]
\draw [dashed] (3,0) node[right]{$b$} -- (0,3) node[right]{$a''$};
\draw [dashed] (-1,2) node[right]{$a$} -- (1,4) node[right]{$a'$};
\fill (3,0) circle (0.05);
\fill (0,3) circle (0.05);
\fill (-1,2) circle (0.05);
\fill (1,4) circle (0.05);
\end{tikzpicture}
\end{center}
\caption{$a'=(a)_b'$. Observers are represented as points; the diagram can be viewed as looking at them ``head on''. Imposing the condition about STL observers guarantees that $a''$ is equally spaced between $a$ and $a'$ (in a Euclidean sense). Essentially, $a'$ is the result of reflecting $a$ over the plane formed by $b$ and $a''$} \label{DualFig}
\end{figure}
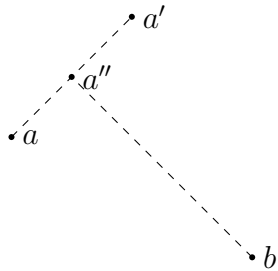

\noindent $a'$ is the \emph{relatable dual} of $a$, with respect to $b$. In corollary \ref{DualCoro}, we will see that $a'$ exists and is unique when $\neg\rho(a,b)$. As noted at the beginning of this section, we do not use the fact that $a'$ is well-defined at any point during the proof of theorem \ref{RepThmFTL} (See figure \ref{DualFig}):
\begin{equation}
\begin{split}
a'=(a)_b'\equiv \neg\rho(a,b)\land a||b\land OP(aa')\land\exists a''[OP(ba'')\land \be_\rho(aa''a') \\
\land \forall \alpha\alpha'\alpha'' cc'c''(l(\alpha\alpha')\land l(\alpha'\alpha'')\land aT\alpha\land a'T\alpha'\land a''T\alpha'' \\ 
\land cT\alpha\land c'T\alpha'\land c''T\alpha'' \land \stl(c)\land c||c'\land c||c'' \Rightarrow \ed_\rho(cc'',c''c'))]
\end{split}
\end{equation}

\noindent $b$ is \emph{between} $a$ and $c$:
\begin{equation}
\be^{FTL}(abc)\equiv\be_\rho(abc)\lor\be_\rho((a)_b'b(c)_b')
\end{equation}

Note that we are using a bit of shorthand here; by this statement we mean:
\begin{equation}
\be^{FTL}(abc)\equiv\be_\rho(abc)\lor\exists a'c'[a'=(a)_b'\land c'=(c)_b' \land\be_\rho((a)_b'b(c)_b')]
\end{equation}

\noindent The distance between $a$ and $b$ is the same as the distance between $c$ and $d$:
\begin{equation}
\ed^{FTL}(ab,cd)\equiv \ed_\rho(ab,cd)\lor\ed_\rho(a(b)_a',c(d)_c')
\end{equation}

\noindent $\alpha$ and $\beta$ are \emph{simultaneous} with respect to $a$ (Figure \ref{SimDefFTL}):
\begin{equation}
\begin{split}
\sm_c^{FTL}(\beta,\beta')\equiv \exists a \gamma\delta[a||c \land aT\gamma\land aT\delta\land \gamma\ne\delta\land l(\beta\gamma) \\ \land l(\beta\delta)\land l(\beta'\gamma) 
\land l(\beta'\delta)]\lor \beta=\beta'
\end{split}
\end{equation}

\begin{figure}
\begin{center}
\begin{tikzpicture}[scale=1]
\draw (0,0) node[right]{$a$} -- (0,6);
\draw (5,0) node[right]{$c$} --  (5,6);
\draw[dashed] (0,1) node[right]{$\gamma$} -- (2,3) node[right] {$\beta'$} 
 -- (0,5) node[right]{$\delta$}
 -- (-2,3) node[left] {$\beta$}
 -- (0,1); 
\fill (2,3) circle (0.05);
\fill (-2,3) circle (0.05);
\fill (0,1) circle (0.05);
\fill (0,5) circle (0.05);
\end{tikzpicture}
\end{center}
\caption{$\sm_c^{FTL}(\beta,\beta')$} \label{SimDefFTL}
\end{figure}

\noindent The difference in time between events $\alpha_0$ and $\alpha_1$ is the same as the difference in time between events $\beta_0$ and $\beta_1$, according to $a$. The only change made from $\Delta$ to $\Delta^{FTL}$ is replacing instances of $\sm$ with $\sm^{FTL}$. (Figure \ref{TimeDifDef}):
\begin{equation}
\begin{split}
\Delta_a^{FTL}(\alpha_0\alpha_1,\beta_0\beta_1)\equiv \exists b \alpha_0'\alpha_1'\beta_0'\beta_1'\{aT\alpha_0'\land aT\alpha_1'\land aT\beta_0'\land aT\beta_1' \\ 
\land \sm_a^{FTL}(\alpha_0,\alpha_0')\land\sm_a^{FTL}(\alpha_1,\alpha_1')\land\sm_a^{FTL}(\beta_0,\beta_0') \\ \land\sm_a^{FTL}(\beta_1,\beta_1')
\land b||a \land [ (\alpha_0'=\alpha_1'\land\beta_0'=\beta_1')\lor(\alpha_0'\ne\alpha_1'\land\beta_0'\ne\beta_1'\land \\
\exists \alpha_2\beta_2(bT\alpha_2\land bT\beta_2\land l(\alpha_0',\alpha_2) \land l(\alpha_1',\alpha_2)\land l(\beta_0',\beta_2)\land l(\beta_1',\beta_2)))]\}
\end{split}
\end{equation}

Finally, we define $\tau^{FTL}$ (Figure \ref{TauDef}). When $c=\tau^{FTL}_b(\alpha_1,\alpha_2)$, we consider the unique observer $a$ that transmits $\alpha_1$ and $\alpha_2$. Then $c$ is the unique observer that transmits some signal $\gamma$ such that $\nd(\gamma)=\alpha_2$ and $\sm_a^{FTL}(\beg(\gamma),\alpha_1)$, and such that $\be^{FTL}(bac)$. As with the previous definition, the only change made from $\tau$ to $\tau^{FTL}$ is replacing $\be$ with $\be^{FTL}$ and $\sm$ with $\sm^{FTL}$. From Lemma \ref{AxLemma2FTL}, we will see that $c$ is uniquely defined when it exists:
\begin{equation}
\begin{split}
c=\tau_b^{FTL}(\alpha_1,\alpha_2)\equiv\exists a \gamma[\stl(a)\land\be^{FTL}(bac)\land aT\alpha_1\land aT\alpha_2 \land \\
cT\gamma\land\sm_a^{FTL}(\beg(\gamma),\alpha_1)\land\nd(\gamma)=\alpha_2\land a\ne b\land\alpha_1\ll\alpha_2]
\end{split}
\end{equation}

\section{Axioms for \textsc{SimpleRelFTL}} \label{axiomsFTL}
In this section we give the axioms of \textsc{SimpleRelFTL} by making some small modifications to \textsc{SimpleRel}. First, we change some of the defined relation symbols that are used. We replace all instances of $\be,\ed,\sm,\Delta$, and $\tau$ with the newly defined $\be^{FTL},$ $\ed^{FTL}$, $\sm^{FTL}$, $\Delta^{FTL}$, and $\tau^{FTL}$. Second, we weaken the statements of \textsc{AxGeo}, \textsc{AxIso}, and \textsc{AxUnSi}. For each of the axioms in \textsc{SimpleRelFTL}, if this axiom is a modified version of an axiom in \textsc{SimpleRel}, we append the superscript $FTL$ to its name. So \textsc{AxSTL} becomes \textsc{AxSTL}\sftl.

Finally, we add the four axioms \textsc{AxFTL1}, \textsc{AxFTL2}, \textsc{AxFTL3}, and \textsc{AxFTL4}. \textsc{AxFTL1} guarantees that there are no lightspeed observers, and \textsc{AxFTL3} guarantees the existence of many FTL observers. Without \textsc{AxFTL3} in particular, any model of \textsc{SimpleRel} would be a model of \textsc{SimpleRelFTL}. There would also be models of \textsc{SimpleRelFTL} that include FTL observers. Including this axiom ensures that all the models of \textsc{SimpleRelFTL} are elementary equivalent, and allows us to give a nice representation theorem for \textsc{SimpleRelFTL}.

\medbreak\textsc{Elementary Geometry Axioms (AxGeo}\sftl): For each axiom that we included in \textsc{AxGeo} in \textsc{SimpleRel}, we add the condition that $a$ is STL, and we replace instances of $\be$ and $\ed$ with $\be^{FTL}$ and $\ed^{FTL}$ . For instance, Axiom 1 of $\mathfrak{E}$ (the first order axiomatization of elementary geometry) becomes:
\begin{equation}
\forall a \{ \stl(a)\Rightarrow \forall xy[ x||a \land y||a \Rightarrow(\be^{FTL}(xyx)\Rightarrow(x=y))]\}
\end{equation}

\textsc{Isotropy Axiom (AxIso}\sftl): We add the condition that $a$ is STL. Given a STL observer $a$ and an event $\beta$, there is a signal from $a$ to $\beta$ and a signal from $\beta$ to $a$.
\begin{equation}
\begin{split}
\forall a\beta[(\stl(a)\land\ev(\beta))\Rightarrow\exists\gamma_1\gamma_2(\beta=\nd(\gamma_1)\land aT\gamma_1 \\ \land\beta=\beg(\gamma_2)\land aR\gamma_2)]
\end{split}
\end{equation}

\textsc{Uniformity~of~Signal~Axiom (AxUnSi}\sftl): Signals travel in straight lines.
\begin{equation}
\begin{split}
\forall abc\alpha\beta\gamma[STL(a)\land\be^{FTL}(abc)\land aT\alpha\land bT\beta \\
\land cT\gamma\land L(\alpha,\gamma)\Rightarrow (L(\alpha,\beta)\Leftrightarrow L(\beta,\gamma))]
\end{split}
\end{equation}

\textsc{Faster-than-light observers axiom 1 (AxFTL1)}: Here we give the first of four assumptions about FTL observers. \textsc{AxFTL1} says that there are no lightspeed observers. 
\begin{equation}
\forall a(\stl(a)\lor\ftl(a))
\end{equation}

\textsc{Faster-than-light observers axiom 2 (AxFTL2)}: If $a$ meets $b$ and $c$, and one of $b,c,d$ is between the other two, then $a$ meets $d$. While this is an assumption we make about all observers, it guarantees that FTL observers don't have any ``gaps'' in their worldlines. This fact follows from \textsc{AxIso} in the proof of theorem \ref{RepThmSTL}. This assumption is needed to guarantee that this holds for FTL observers, since in \textsc{SimpleRelFTL}, \textsc{AxIso}\sftl does not guarantee anything about FTL observers.
\begin{equation}
\begin{split}
\forall abcd \beta \gamma_1\gamma_2\ [(\be^{FTL}(bcd)\lor\be^{FTL}(cbd) \lor\be^{FTL}(bdc)) \\
\land M(ab)\land M(ac) 
\Rightarrow M(ad))]
\end{split}
\end{equation}

\textsc{Faster-than-light observers axiom 3 (AxFTL3)}: Given an observer $a$ and signal $\beta$ such that $aR\beta$, for any event $\gamma_2$ transmitted by $a$ that occurs \emph{before} $\nd(\beta)$, there is an observer $b$ that transmits $\beta$ and $\gamma_2$. Note the similarity of this assumption to \textsc{AxLim}. 
\begin{equation}
\begin{split}
\forall a \beta \gamma_1\gamma_2 [(\gamma_1=\nd(\beta)  \land\gamma_2\ll\gamma_1  \land aT\gamma_1\land aT\gamma_2 ) \\ \Rightarrow\exists b (bT\beta\land bT\gamma_2)] 
\end{split}
\end{equation}

\textsc{Faster-than-light observers axiom 4 (AxFTL4)}: All FTL observers transmit at least two signals. This fact was guaranteed by \textsc{AxGeo} in \textsc{SimpleRel}. We include this assumption now because nothing was assumed about the geometry of the world from a FTL observer's perspective.
\begin{equation}
\begin{split}
\forall a (FTL(a)\Rightarrow \exists \alpha_1\alpha_2(aT\alpha_1\land aT\alpha_2\land \beg(\alpha_1)\ne\beg(\alpha_2)))
\end{split}
\end{equation}

\section{Useful Lemmas for \textsc{SimpleRelFTL}} \label{FTLLemmas}

In this section, we prove several lemmas which will be used in the representation theorem of \textsc{SimpleRelFTL}. Most importantly, we show that if $a$ is STL, then the definitions of betweenness, equidistance, and simultaneity given in this chapter are equivalent to those given in the previous one. We also note that analogues of many of the lemmas proved in section \ref{STLLemmas} hold.

\begin{lemma} \label{ParLemmaFTL} Assume \textsc{AxGeo}\sftl and \textsc{AxSTL}\sftl. Then the relation $||$ is an equivalence relation between STL observers. \end{lemma}
\proof See the proof of lemma \ref{ParLemma}. \proofbox

\begin{coro} \label{AxLemma1FTL} Assume \textsc{AxRR}\sftl and \textsc{AxGeo}\sftl. If $a$ is an STL observer and $\gamma$ is an event, there exists a \emph{unique} observer $b$ such that $bT\gamma$ and $b||a$.
\end{coro}
\proof See the proof of corollary \ref{AxLemma1}. \proofbox

\begin{lemma} \label{DefEquivLemma} Assume \textsc{AxIso}\sftl, \textsc{AxSTL}\sftl, \textsc{AxGeo}\sftl, \textsc{AxEv}\sftl, \\ \textsc{AxUnSi}\sftl,
\textsc{AxTime}\sftl, and \textsc{AxLim}\sftl. If $a$ is STL, then for all observers $b,c,d$ and events $\alpha,\beta$, $\be^{FTL}(abc)\Leftrightarrow\be(abc)$, $\ed(ab,cd)^{FTL}\Leftrightarrow\ed(ab,cd)$, $\sm_a^{FTL}(\alpha,\beta)\Leftrightarrow\sm_a(\alpha,\beta)$. What is more, for any events $\alpha_0,$ $\alpha_1,$ $\beta_0,$ $\beta_1$, 
\begin{equation}
\Delta_a(\alpha_0\alpha_1,\beta_0\beta_1)\Leftrightarrow\Delta_a^{FTL}(\alpha_0\alpha_1,\beta_0\beta_1),
\end{equation}
 and for any observers $b,c$, 
 \begin{equation}
 c=\tau_b(\alpha_1,\alpha_2)\Leftrightarrow c=\tau_b^{FTL}(\alpha_1,\alpha_2).
 \end{equation}
\end{lemma}
\proof Because $a$ is STL, \textsc{AxIso}\sftl ~holds for $a$, and so for any observer $b$, $\rho(a,b)$. To prove the first two statements, it then suffices to show that $\be_\rho(abc)\Leftrightarrow\be(abc)$ and $\ed_\rho(ab,cd)\Leftrightarrow\ed(ab,cd)$. Because of how $\Delta^{FTL}$ and $\tau^{FTL}$ are defined, to prove the last two statements it suffices to prove the first and the third.

Assume that $\be_\rho(abc)$. By \textsc{AxSTL}\sftl and lemma \ref{ParLemmaFTL}, $a||b||c$. By the definition of $\be_\rho$ and $\be$, either $\be(abc)$ or $\be(cba)$. And if $\be(cba)$, there exists events $\alpha,\beta,\gamma$ such that $aT\alpha$, $bT\beta$, $cT\gamma$, $L(\gamma,\beta)$, $L(\gamma,\alpha)$, and $L(\beta,\alpha)$. Since $a$ is STL, by \textsc{AxSTL}\sftl ~$b$ and $c$ are STL. By \textsc{AxIso}\sftl and \textsc{AxEv}\sftl there exists events $\gamma'$ and $\alpha'$ such that $cT\gamma'$, $aT\alpha'$, $L(\beta,\gamma')$, and $L(\alpha',\gamma')$. Since $\be_\rho(abc)$ and hence $\be^{FTL}(abc)$, we may apply \textsc{AxUnSi}\sftl \\~to find that $L(\alpha',\beta)$ and hence $\be(abc)$. Similarly, noting the definitions of $\be$ and $\be_\rho$, after applying lemma \ref{ParLemmaFTL} we have that $\be(abc)\Rightarrow\be_\rho(abc)$. Thus, $\be^{FTL}(abc)\Leftrightarrow\be(abc)$.

Now, assume that $\ed_\rho(ab,cd)$. If $a=b$ and $c=d$ the proof is complete, because by \textsc{AxGeo}\sftl there exists an observer $e$ such that $\ed_\rho(ea,ed)$ and hence there exist signals that satisfy the definition of $\ed(ab,cd)$. And because $a$ is STL, by \textsc{AxTime}\sftl and \textsc{AxLim}\sftl, $\neg OP(ab)$ and $\neg OP(cd)$.
 
We may assume then that there exists signals $\alpha_1\ne\alpha_2,\beta,\gamma_1\ne\gamma_2,$ and $\delta$ such that $l(\alpha_1\beta)$, $l(\alpha_2\beta)$, $l(\gamma_1\delta)$, $l(\gamma_2\delta)$, $aT\alpha_1$, $aT\alpha_2$, $bT\beta$, $cT\gamma_1$, $cT\gamma_2$, and $dT\delta$, and that there exists an observer $e$ as in the definition. We may also assume that for every signal $\epsilon$, if $e$ transmits $\epsilon$ then $l(\epsilon\alpha_1)\Leftrightarrow l(\epsilon\gamma_1)$ and $l(\epsilon\alpha_2)\Leftrightarrow l(\epsilon\gamma_2)$. By \textsc{AxSTL}\sftl, $b,c,d,e$ are STL and from lemma \ref{ParLemmaFTL} we have that $a||b||c||d||e$. \textsc{AxIso}\sftl then gives us that there exists signals $\alpha_3,\alpha_4$ from $\alpha_1,\alpha_2$ to $e$. By \textsc{AxEv}\sftl and because $\ed_\rho(ab,cd)$, $l(\nd(\alpha_3),\gamma_1)$. Thus, $L(\nd(\alpha_3),\gamma_1)$ or $L(\gamma_1,\nd(\alpha_3))$. But if $L(\nd(\alpha_3),\gamma_1)$, because $e$ is STL there exists some signal $\gamma_3$ such that $\beg(\gamma_3)=\gamma_1$ and $eR\gamma_3$. Then $\nd(\alpha_3)\ll\nd(\gamma_3)$ and so by \textsc{AxTime}\sftl, $\nd(\alpha_3)\ne\nd(\gamma_3)$. Now, because $\ed_\rho(ab,cd)$, $L(\alpha_1,\nd(\gamma_3))$ or $L(\nd(\gamma_3),\alpha_1)$. If $L(\alpha_1,\nd(\gamma_3))$, then there are two distinct signals from $\alpha_1$ to $e$, contradicting the fact that $c$ is STL. And if $L(\nd(\gamma_3),\alpha_1)$, then $\nd(\gamma_3)\ll\nd(\alpha_3)$, contradicting \textsc{AxTime}\sftl. Thus, $L(\gamma_1,\nd(\alpha_3))$. A similar argument shows that $L(\gamma_2,\nd(\alpha_4))$. Now, by \textsc{AxIso}\sftl there exists events $\epsilon_1,\epsilon_2$ transmitted by $e$ such that $L(\epsilon_1\alpha_1)$ and $L(\epsilon_2\alpha_2)$. We then also have that $l(\epsilon_1\gamma_1)$ and $l(\epsilon_2\gamma_2)$. Since $c$ is STL, it follows that $L(\epsilon_1\gamma_1)$ and $L(\epsilon_2\gamma_2)$. Letting $\nd(\alpha_3)=\epsilon_1'$ and $\nd(\alpha_4)=\epsilon_2'$, we find that the conditions for $\ed(ab,cd)$ are satisfied.

Assume that $\ed(ab,cd)$, and let $e,$$\alpha_1,$$\alpha_2,$$\beta,$$\gamma_1,$$\gamma_2,$$\delta,$$\epsilon_1,$$\epsilon_1',$$\epsilon_2,$ and $\epsilon_2'$ be as they are in the definition of $\ed$. Then by \textsc{AxSTL}\sftl, $b,c,d,e$ are STL and by lemma \ref{ParLemmaFTL} $a||b||c||d||e$. If $a=b$ and $c=d$ then the proof is complete. Note $a\ne b$ iff $c\ne d$, and so when $a\ne b$, $\alpha_1\ne\alpha_2$ and $\gamma_1\ne\gamma_2$. Since $a,c,$ and $e$ are STL, $\epsilon_1,$$\epsilon_1',$$\epsilon_2,$ and $\epsilon_2'$ give all the events transmitted by $e$ that are lightlike with $\alpha_1,$$\alpha_2,$$\gamma_1,$ and $\gamma_2,$ and so the definition for $\ed_\rho(ab,cd)$ is satisfied. Thus, if $a$ is STL, $\ed(ab,cd)\Leftrightarrow\ed^{FTL}(ab,cd)$. 

Finally, (letting now $c$ be STL and $\beta,\beta'$ be two events), we show that $\sm_c(\beta,\beta')\Leftrightarrow\sm_c^{FTL}(\beta,\beta')$. Assume that $\sm_c(\beta,\beta')$ and let $a,\gamma_1,\gamma_2,\gamma_3,\gamma_4$ be as in the definition of $\sm$. If $\beta=\beta'$ the proof is complete. Otherwise, since $a$ is STL by \textsc{AxSTL}\sftl, by \textsc{AxTime}\sftl it follows that $\beg(\gamma_1)\ne\nd(\gamma_3)$. Setting $\beg(\gamma_1)=\gamma$ and $\nd(\gamma_3)=\delta$ gives that $\sm_c^{FTL}(\beta,\beta')$.

Assume that $\sm_c^{FTL}(\beta,\beta')$. If $\beta=\beta'$, $\sm_c(\beta,\beta')$ because $c$ is STL. Otherwise, note that, since $\be^{FTL}(aa'a'')\Leftrightarrow\be(aa'a'')$ for all $a,a',a''$ such that $a||a'||a''$ and $\stl(a)$, the result of lemma \ref{STLLemma1} holds for all $a||c$ (see lemma \ref{STLLemma1FTL}, below). Taking $a,\gamma,\delta$ as in the definition of $\sm_c^{FTL}$, we then have that either $L(\gamma\beta)$, $L(\beta\delta)$, $L(\gamma\beta')$, and $L(\beta'\delta)$, or $L(\gamma\beta)$, $L(\beta\delta)$, $L(\beta'\gamma)$, and $L(\delta\beta')$. In the first case, $\sm_c(\beta,\beta')$ holds. In the second case, $\gamma \ll \delta$ and $\delta \ll \gamma$, contradicting \textsc{AxTime}\sftl. Thus, $\sm_c^{FTL}(\beta,\beta')\Rightarrow\sm_c(\beta,\beta')$.\proofbox

\begin{lemma} \label{UniquenessOfSignalsLemmaFTL} Assume \textsc{AxEv}\sftl, \textsc{AxSTL}\sftl, and \textsc{AxRR}\sftl. Then if two signals $\alpha$ and $\beta$ share endpoints, then $\alpha=\beta$. \end{lemma}

\proof See the proof of lemma \ref{UniquenessOfSignalsLemma}. \proofbox

\begin{lemma} \label{TwoEventsLemmaFTL} Assume \textsc{AxIso}\sftl, \textsc{AxSTL}\sftl, \textsc{AxGeo}\sftl, \textsc{AxUnSi}\sftl, \newline 
\textsc{AxTime}\sftl, \textsc{AxLim}\sftl, \textsc{AxEv}\sftl, and \textsc{AxFTL4}. If $a$ is an observer, there are at least two distinct events $\alpha_1$ and $\alpha_2$ that are transmitted by $a$. \end{lemma}

\proof If $a$ is STL, then by applying lemma \ref{DefEquivLemma} the proof of lemma \ref{TwoEventsLemma} holds for $a$. If $a$ is not STL, then the statement holds by \textsc{AxFTL4} and \textsc{AxEv}\sftl. \proofbox

\begin{coro} \label{ObserversEventsLemmaFTL} Assume \textsc{AxIso}\sftl, \textsc{AxSTL}\sftl, \textsc{AxGeo}\sftl, \textsc{AxUnSi}\sftl, \textsc{AxTime}\sftl, \textsc{AxLim}\sftl, \textsc{AxEv}\sftl, \textsc{AxFTL4}, \textsc{AxRR}\sftl, and \textsc{AxOb!}\sftl. Then two observers $a$ and $b$ are identical iff for every signal $\alpha$, $aT\alpha\Leftrightarrow bT\alpha$. \end{coro}

\proof If $a=b$ then they transmit all the same signals. For observers $a$ and $b$, assume that given a signal $\alpha$, $aT\alpha\Leftrightarrow bT\alpha$. From Lemma \ref{TwoEventsLemmaFTL}, $a$ transmits at least two events. Thus, by \textsc{AxOb!}\sftl, $a=b$. \proofbox

\begin{lemma} \label{STLLemma1FTL} Assume \textsc{AxIso}\sftl, \textsc{AxSTL}\sftl, \textsc{AxGeo}\sftl, \textsc{AxUnSi}\sftl, \newline 
\textsc{AxTime}\sftl, \textsc{AxLim}\sftl, \textsc{AxEv}\sftl, and \textsc{AxRR}\sftl. Then if $a$ is STL, then for all events $\gamma$, there exists a unique signal $\beta$ such that $aT\beta$ and $\nd(\beta)=\gamma$ (Figure \ref{STLLemma1Fig}).
\end{lemma}

\proof Applying lemma \ref{DefEquivLemma}, see the proof of lemma \ref{STLLemma1}. \proofbox

\begin{lemma} \label{STLEventsLemmaFTL} Assume \textsc{AxIso}\sftl, \textsc{AxSTL}\sftl, \textsc{AxGeo}\sftl, \textsc{AxUnSi}\sftl, \newline 
\textsc{AxTime}\sftl, \textsc{AxLim}\sftl, \textsc{AxEv}\sftl, and \textsc{AxRR}\sftl. If $a$ is STL, $aT\alpha_1$, $aT\alpha_2$, and $\sm_a(\alpha_1,\alpha_2)$, then $\alpha_1=\alpha_2$. \end{lemma}

\proof Applying lemma \ref{DefEquivLemma}, see the proof of lemma \ref{STLEventsLemma}. \proofbox

\begin{lemma} \label{AxLemma2FTL} Assume \textsc{AxIso}\sftl, \textsc{AxSTL}\sftl, \textsc{AxGeo}\sftl, \textsc{AxUnSi}\sftl, \newline
\textsc{AxTime}\sftl, \textsc{AxLim}\sftl, \textsc{AxEv}\sftl, \textsc{AxOb!}\sftl, and \textsc{AxPoInd}\sftl. Then if $b=\tau_c(\alpha_1,\alpha_2)$ and $b'=\tau_c(\alpha_1,\alpha_2)$, then $b=b'$.
\end{lemma}

\proof Applying lemma \ref{DefEquivLemma}, see the proof of lemma \ref{AxLemma2}. \proofbox

\begin{lemma} \label{SimEventsLemmaFTL} Assume \textsc{AxIso}\sftl, \textsc{AxSTL}\sftl, \textsc{AxGeo}\sftl, \textsc{AxUnSi}\sftl, \newline 
\textsc{AxTime}\sftl, \textsc{AxLim}\sftl, \textsc{AxEv}\sftl,  \textsc{AxOb!}\sftl, \textsc{AxPoInd}\sftl, and \newline
\textsc{AxTiInd}\sftl. For all events $\beta$ and all slower than light observers $a$ there exists a unique event $\alpha$ such that $aT\alpha$ and $\sm_a(\alpha,\beta)$ (Figure \ref{SimEventsLemmaFig}).
\end{lemma}

\proof Applying lemma \ref{DefEquivLemma}, see the proof of lemma \ref{SimEventsLemma}. \proofbox

\section{Representation theorem for \textsc{SimpleRelFTL}}
\label{RepThmFTL}
In this section, we prove a representation theorem for \textsc{SimpleRelFTL}. Most of it relies on the proof of the representation theorem for \textsc{SimpleRel}. We first define the class of structures that will be our models for \textsc{SimpleRelFTL}

As we did in section \ref{RCOFmodels}, we let $\mathbf{F}$ be a real-closed ordered field, $\lambda$ be the Minkowski norm. We say that a line $\mathbf{a\subset F^4}$ is a \emph{spacelike line} if there exists $\mathbf{v,w\in F^4}$ such that $\lambda(\mathbf{w})\mathbf{>0}$ and $\mathbf{a=\{v+tw:t\in F\}}$. We define $\M_F^{FTL}=\mathbf{\langle Ob_F^{FTL}, Si_F^{FTL},T_F^{FTL},R_F^{FTL}\rangle}$, where $\mathbf{Si_F^{FTL}=Si_F}$, $\mathbf{T_F^{FTL}=T_F}$, $\mathbf{R_F^{FTL}=R_F}$, and $\mathbf{Ob_F^{FTL}\subset\pow(F^4)}$ be the set of all timelike \emph{and spacelike} lines of $\mathbf{F^4}$.

\begin{theo} $\M^{FTL}=\mathbf{\langle Ob^{FTL}, Si^{FTL},T^{FTL},R^{FTL}\rangle}$ is a model of \textsc{SimpleRelFTL} iff there is a real-closed ordered field $\mathbf{F}$ such that $\M^{FTL}\cong\M_F^{FTL}$. \end{theo}

\proof
$\Leftarrow$: One checks that $\M_\mathcal{R}^{FTL}$ satisfies each of the axioms in \textsc{SimpleRelFTL}. Since $\mathbf{F}$ is elementary equivalent to $\mathcal{R}$, it follows that $\M_F^{FTL}$ is elementary equivalent to $\M_\mathcal{R}^{FTL}$, and hence $\M_F^{FTL}$ is a model of \textsc{SimpleRelFTL}.

$\Rightarrow$: Applying lemma \ref{DefEquivLemma}, we may follow the proof of theorem \ref{RepThmSTL}, up to the point where we have defined $\psi_{Ob},\psi_{Si}$, shown that $\psi_{Si}:\mathbf{Si^{FTL}\to Si_F^{FTL}}$ is a bijection and that $\psi_{Ob}:\M^{FTL}\to\pow(\mathbf{F^4})$ maps observers to subsets of straight lines in $\mathbf{F^4}$. And if $\mathbf{b}$ is STL, then by \textsc{AxIso}\sftl and the fact that $\psi_{Si}$ is a bijection, $\psi_{Ob}(\mathbf{b})$ is a timelike line in $\mathbf{F^4}$. 

If $\mathbf{b}$ is not STL, then by \textsc{AxFTL1}, $\mathbf{b}$ is FTL. In this case, because of the structure of $\mathbf{F^4}$ and because $\psi_{Si}$ is a bijection, $\psi_{Ob}(\mathbf{b})$ must lie within a spacelike line $\mathbf{l(b)}\subset\mathbf{F^4}$. By lemma \ref{TwoEventsLemmaFTL}, $\mathbf{b}$ transmits at least two events. Consider the observers $\mathbf{c,c'\in P_a}$ that transmit these events; $\mathbf{c,c'}$ exist by \textsc{AxRR}\sftl. Now, given an event $\mathbf{\alpha}$ such that $\psi(\mathbf{\alpha})\in \mathbf{l(b)}$, there exists an observer $\mathbf{c''\in P_a}$ such that $\mathbf{c''T \alpha}$. Also, one of $\psi_{Ob}(\mathbf{c}),\psi_{Ob}(\mathbf{c'}),\psi_{Ob}(\mathbf{c''})$ is between the other two, since they are all vertical lines in $\mathbf{F^4}$ that meet $\mathbf{l(b)}$. And since $\psi_{Si}$ is a bijection, we find that one of $\mathbf{c,c',c''}$ is between the other two. Hence, by \textsc{AxFTL2} $\mathbf{M(b,c'')}$ and so $\mathbf{bT\alpha}$. Thus, $\psi(\mathbf{b})\in\mathbf{Ob_F^{FTL}}$ and ${\rm{RANGE}}(\psi_{Ob})=\mathbf{Ob_F^{STL}}$. And $\psi_{Ob}$ is one-to-one because $\psi$ is one-to-one.

Finally, we show that $\psi_{Ob}$ is onto. That is, if $\mathbf{a\in Ob_F^{FTL}}$ then there exists an observer $\mathbf{a'\in Ob^{FTL}}$ such that $\psi_{Ob}(\mathbf{a'})=\mathbf{a}$. If $\mathbf{a}$ is STL, then the proof follows the one given in theorem \ref{RepThmSTL}. If $\mathbf{a}$ is FTL, then the proof also follows this one, except instead of using \textsc{AxLim}\sftl, we use \textsc{AxFTL3}. Thus, $\psi_{Ob}$ is onto and $\M^{FTL}\cong\M_F^{FTL}$. \proofbox

\begin{coro} \label{DualCoro} It is provable from \textsc{SimpleRelFTL} that the relatable dual of $a$ with respect to $b$, $(a)_b'$ exists and is uniquely defined when $a||b$ and $\neg\rho(a,b)$. \end{coro}
\proof Since this statement holds true in $\M_{\mathcal{R}}^{FTL}$, by the representation theorem for \textsc{SimpleRelFTL} it holds in every model of \textsc{SimpleRelFTL}. By the Completeness Theorem (Theorem \ref{ComTheo}), the statement is provable from \textsc{SimpleRelFTL}. \proofbox

\begin{coro} The set of consequences of \textsc{SimpleRelFTL} is complete and decidable. \end{coro}
\proof See the proof of corollary \ref{Coro2STL} \proofbox

\begin{note} \end{note} In section \ref{definitions}, we gave a number of definitions of geometric and physical predicates. We then gave the axiom system \textsc{SimpleRel} using these definitions, and proved a representation theorem about \textsc{SimpleRel}. However, we noted that some of these definitions were only meaningful for STL observers. In the section \ref{definitionsFTL}, we gave definitions of these terms that are meant to be meaningful for both STL and FTL observers. Let \textsc{SimpleRel}$^*$ contain the same axioms as \textsc{SimpleRel}, except that we replace instances of $\be,\ed,\sm,\Delta,$ and $\tau$ with $\be^{FTL},\ed^{FTL},\sm^{FTL},\Delta^{FTL},$ and $\tau^{FTL}$. By slightly modifying the proofs of the lemmas and the theorem in this chapter, one may show that \textsc{SimpleRel}$\models$\textsc{SimpleRel}$^*$ and \textsc{SimpleRel}$^*\models$\textsc{SimpleRel}. The definitions of the predicates $\be^{FTL},\ed^{FTL},\sm^{FTL},\Delta^{FTL},$ and $\tau^{FTL}$ are then generalizations of the definitions of the corresponding predicates used in \textsc{SimpleRel}, and we could have equally given \textsc{SimpleRel} using these new definitions.


\end{document}